\documentclass{article}
\usepackage{amssymb,latexsym,amsmath}
\usepackage{graphicx}
\newtheorem{theorem}{Theorem}[section]

\newtheorem{proposition}[theorem]{Proposition}
\newtheorem{corollary}[theorem]{Corollary}
\newtheorem{definition}[theorem]{Definition}
\newcommand{\supp}{\rm supp}

\title{Asymptotics and zeros of Sobolev orthogonal polynomials on unbounded supports}

\author{Francisco Marcell{\'a}n\thanks{Research partially supported by Direcci\'{o}n General de
Investigaci\'{o}n (Ministerio de Ciencia y Tecnolog\'{\i}a) of
Spain under grant grant BFM2003--06335--C03--02 and INTAS Research
Network NeCCA INTAS 03--51--6637.}
\\ Departamento de Matem\'aticas \\
Universidad Carlos III de Madrid, Spain \and  Juan Jos{\'e} Moreno
Balc{\'a}zar\thanks{Research partially supported by Ministerio de
Educaci{\'o}n y Ciencia of Spain under grant MTM2005-08648-C02-01 and
Junta de Andaluc\'{\i}a (FQM 229 and FQM 481).} \\ Departamento de
Estad{\'\i}stica y Matem{\'a}tica Aplicada
\\Universidad de Almer{\'\i}a, Spain
 \\ and \\ Instituto Carlos I de F{\'\i}sica
Te{\'o}rica y Computacional \\ Universidad de Granada, Spain }

 \date{}

\begin{document}

\maketitle

\begin{abstract}
In this paper we present a survey about analytic properties of
polynomials orthogonal with respect to a weighted Sobolev inner
product such that the vector of measures has an unbounded support.
In particular, we focus on the asymptotic behaviour of such
polynomials as well as in the distribution of their zeros. Some
open problems as well as some  directions for  future research are
formulated.

\end{abstract}

\medskip

\noindent {\it MSC:} Primary 42C05; Secondary 33C45

\medskip

\noindent {\it Key words}: Sobolev orthogonal polynomials,
asymptotics, zeros.

\section{Introduction}

The study of polynomials orthogonal with respect to a Sobolev
inner product
\begin{equation} \label{pr-sob-nder}
(p,q)=\sum_{k=0}^N\int p^{(k)}(x)\, q^{(k)}(x)\, d\mu_k(x)
\end{equation}
where $\mu=(\mu_0, \mu_1, \ldots, \mu_N)$ is a vector of positive
Borel measures supported on the real line and $p,q\in \mathbb{P},$
where $\mathbb{P}$ denotes the set of polynomials with real
coefficients, appears in many areas of Applied Mathematics. In
smooth data fitting (\cite{lew}) they constitute the natural basis
for the approximating subspaces. In Fourier expansions, the
analysis of Gibbs phenomena can be performed in a successful way
using such polynomials as the numerical experiments given by A.
Iserles et al. \cite{ise-con} show. Finally, the analysis of
spectral methods for boundary value problems based on such
polynomials is more competitive than the consideration of standard
orthogonal polynomials.

The subject has attracted increased attention in the last 15
years. As a sample of it, the extensive bibliography by F.
Marcell{\'a}n and A. Ronveaux \cite{mar-ron} shows the different
topics of interest for researchers. Most of the contributions
therein are related to vector of measures $\mu$ such that their
components are classical or are like--classical measures, i.e,
Jacobi, Laguerre, and Hermite measures. This is a natural
consequence of the well known behaviour of the corresponding
sequences of orthogonal polynomials from an analytic point of
view. In particular,  attention was focused on $N=1$ in order to
understand the simplest cases.

When an inner product like (\ref{pr-sob-nder}) is considered, two
basic properties of standard orthogonal polynomials ($N=0$) are
lost. First, the polynomials orthogonal with respect to
(\ref{pr-sob-nder}) with $N>0$ do not satisfy a three--term
recurrence relation that constitutes a useful tool in the analysis
of the ratio asymptotics of standard orthogonal polynomials.
Second, the zeros of polynomials orthogonal with respect to
(\ref{pr-sob-nder}) with $N>0$ are not, in general, real and
simple as in the standard case. They can be complex and even, if
real, they can be located outside the convex hull of the support
of the vector of measures $\mu.$

When $\mu$ has  bounded support, the asymptotic behaviour of the
polynomials orthogonal with respect to  (\ref{pr-sob-nder}) has
been extensively studied by A. Mart{\'\i}nez--Finkelshtein, G.
L{\'o}pez--Lagomasino among others. In particular, relative
asymptotics, $n$th root asymptotics, and strong asymptotics have
been deeply analyzed.

Nevertheless, when $\mu$ is a vector of measures with unbounded
support, despite the examples related to Laguerre weights (the so
called coherent pairs) or Hermite weights (the so called
symmetrically coherent pairs) very few examples are known. In
\cite{cac-mar-mb} the case of the vector of measures
$d\mu=(e^{-x^4}dx, \lambda e^{-x^4}dx )$ is analyzed. In
particular, relative asymptotics of this Sobolev orthogonal
polynomials in terms of the Freud polynomials associated with the
weight function $w(x)=\exp(-x^4)$ are deduced. In a more general
context, in \cite{ger-lub-mar} different kinds of asymptotic
behaviour of Sobolev orthogonal polynomials for exponential
weights are analyzed in the case $N=1.$ In particular, outer and
inner asymptotics in terms of $L_2$ and $L_{\infty}$ norms as well
as strong asymptotics are considered. The key role is played in
these problems by the polynomials orthogonal with respect to the
measure $\mu_1.$

The aim of our paper is to present the  state of the art in the
area of analytic properties of Sobolev polynomials orthogonal with
respect to vector of measures of exponential type supported on
unbounded sets of the real line. We will focus our attention in
the asymptotic behaviour of such polynomials as well as in the
distribution of their zeros, emphasizing the comparison with the
same properties for standard orthogonal polynomials with respect
to the measures $\mu_0$ and $\mu_1.$ Thus we give a survey
following the approach of F. Marcell{\'a}n, M. Alfaro, and M.L. Rezola
(\cite{mar-alf-rez}), H. G. Meijer (\cite{meijer-sur}) and A.
Mart{\'\i}nez--Finkelshtein (\cite{and-sevilla}--\cite{and-survey}).

We hope it will be useful for people interested in this field not
only from a theoretical point of view but taking into account the
potential applications.

The structure of the paper is the following. In Section 2, we deal
with Sobolev inner products (\ref{pr-sob-nder}) for $N=1$ when one
of the measures is Laguerre  (resp. Hermite) and the other one is
its coherent (resp.  symmetrically coherent) companion in the
sense of Meijer classification \cite{meijer}. In the Laguerre case
we also consider a non--diagonal inner product connected with the
framework of the coherence.  In both cases, relative asymptotics
in terms of the corresponding classical polynomials as well as
strong asymptotics on compact subsets of $(0, \infty)$ (resp. real
line)  in Laguerre case (resp. in  Hermite case) are studied. The
distribution of the zeros of such polynomials as well as the
location in terms of the zeros of the corresponding classical
orthogonal polynomials as well as interlacing and separation
properties are deduced. Finally, the behaviour of the zeros of
these Sobolev orthogonal polynomials in terms of the zeros of
Bessel functions is deduced by means of generalized Mehler--Heine
type formulas for Laguerre and Hermite polynomials. Furthermore,
in Proposition \ref{pro-imp} we establish that the natural family
of  polynomials  to compare with these Sobolev polynomials is the
family of polynomials defined as $R_n(x)=\lim_{\lambda \to \infty}
Q_n(x).$ This proves the important role of the parameter
$\lambda.$

In Section 3, we leave the framework of coherence. Here,
asymptotics properties for polynomials orthogonal with respect to
a vector of measures with $d\mu_0=\psi(x)\, \exp(-Q(x))dx, $ $\,
d\mu_1=\lambda\, \exp(-Q(x))dx$ where $Q$ is an even and
continuous function on the real line satisfying some extra
conditions and $\psi$ is a measurable and positive function on a
set of positive measure on $\mathbb{R}$ are studied. Again the
parameter $\lambda$ plays an important role in the asymptotic
behaviour of these Sobolev orthogonal polynomials. This is one of
the main differences with the case of measures of bound support.

Finally, in Section 4 some open problems are posed as well as some
new directions for a further research in the subject are
presented.

\section{The ``classical'' cases}

We consider the Sobolev inner product
\begin{equation} \label{pr-sob-gen}
(p,q)=\int p(x)q(x)d\mu_0(x)+\int p'(x)q'(x)d\mu_1(x)
\end{equation}
where $\mu_i(x)$ are measures supported on subsets $I_i\subseteq
\mathbb{R}.$

We speak of  the ``classical" cases when one of the measures
involved in the inner product (\ref{pr-sob-gen}) is a classical
measure with unbounded support, i.e., the measures corresponding
to either the Laguerre or Hermite weight functions. Indeed, this
situation is the only one  studied for measures with unbounded
supports until the work \cite{cac-mar-mb} for an example of Freud
weights. Three contributions concerning Laguerre and Hermite
weights were given in the 90's (see \cite{mar-mei-per-pin},
\cite{pan}, \cite{mei-per-pin}). These papers are in the framework
of the concept of coherent and symmetrically coherent pairs of
measures introduced in \cite{iserles} being one of the measures
the corresponding to the Laguerre weight. In fact,  in 1997 H.G.
Meijer established the complete classification both of the
coherent and symmetrically coherent pairs \cite{meijer}. In the
case of measures with unbounded support we get:

\bigskip

\noindent \underline{Laguerre weights:}

\bigskip

Case I:

\begin{itemize}

\item{I.1.} \begin{equation} d\mu_0=(x+a)x^{\alpha-1} e^{-x}dx,\,
d\mu_1=x^{\alpha} e^{-x}dx, \quad a\ge 0, \, \, \alpha
>0.\label{ltype11}
\end{equation}

\item{I.2}. \begin{equation} d\mu_0= e^{-x}dx+M\delta (x),\,
d\mu_1=e^{-x}dx,\quad M\ge 0. \label{ltype12} \end{equation}

\end{itemize}

Case II:

 \begin{equation} d\mu_0=x^{\alpha} e^{-x}dx,\, d\mu_1=\frac{x^{\alpha+1}
e^{-x}}{x+a} dx+ M\delta (x+a),\quad a\ge 0, \, \, \alpha
>-1, \quad M\ge 0.\label{ltype2} \end{equation}

\bigskip

\noindent \underline{Hermite weights:}

\bigskip

Case I:
\begin{equation}
d\mu_0=(x^2+a^2) e^{-x^2}dx,\, d\mu_1= e^{-x^2}dx,\quad a\in
\mathbb{R}.\label{htype1} \end{equation}

Case II:
\begin{equation}
 d\mu_0=e^{-x^2}dx,\, d\mu_1=\frac{e^{-x^2}}{x^2+a^2}dx,\quad a\in \mathbb{R}\setminus \{0\}.
\label{htype2} \end{equation}

\bigskip

Notice that the cases II are rational modifications of the
classical measure while the cases I.1 and I for the Laguerre and
Hermite weights, respectively, are polynomial modifications of
these measures. Therefore, it is natural to think that we will
have more problems to establish analytic properties for the
Sobolev orthogonal polynomials of Type II because we will need to
study families of standard orthogonal polynomials with respect to
a rational modification of a measure.

\subsection{Laguerre weights}

The paper \cite{mar-mei-per-pin} is a very particular case when
both measures correspond to the Laguerre weight function. The
authors obtain the asymptotic behaviour of these Sobolev
orthogonal polynomials uniformly on compact subsets of
$\mathbb{C}\setminus [0,\infty).$ In the paper \cite{pan} the
$nth$--root asymptotic of the Laguerre--Sobolev polynomials is
obtained in $\mathbb{C}\setminus [0,2]$ (although along the paper
the sentence ``locally uniformly on compact subsets of $[0,2]$"
appears  as a misprint). Later on,  in \cite{mei-per-pin} the
authors deal with the analytic properties of Laguerre--Sobolev
orthogonal polynomials of both types. They obtain the asymptotics
in $\mathbb{C}\setminus [0, \infty)$. We will introduce these
results in the following items.

\subsubsection{Non--Diagonal Laguerre--Sobolev orthogonal polynomials and  their connection
 with Laguerre--Sobolev of case I}

 In \cite{mar-mb}  the inner product
 \begin{equation}
(p,q)_S=\int_{0}^{\infty} (p,p^{\prime}) \left(
\begin{array}{ll}
1 & \delta \\
\delta & \lambda
\end{array}
\right) \left(
\begin{array}{l}
q \\
q^{\prime}
\end{array}
\right) x^{\alpha }e^{-x}dx \, ,  \label{iplsnd}
\end{equation}
with $\alpha >-1\, $ is considered. There, the assumption
$\lambda, \delta \in \mathbb{R}$ with $\lambda -\delta^2
>0\, $ yields the existence of a unique sequence $\{ Q_n^{(ND)} \}
$ of polynomials orthogonal with respect to (\ref{iplsnd}) up to
normalization. In \cite{mar-mb} the polynomials $Q_n^{(ND)}$ were
considered as monic. Here, we denote by
$L_n^{(\alpha)}(x)=((-1)^n/n!)\, x^n +\ldots $ the classical
Laguerre polynomials and  we define the polynomials $Q_n^{(ND)}$
in such a way they have the same leading coefficient.

This inner product and  the corresponding Sobolev orthogonal
polynomials are related to Laguerre--Sobolev polynomials of case
I.1 whose measures are given by (\ref{ltype11}) through a
straightforward integration by parts. Indeed, for $\alpha>0,$ we
have
$$ (p,q)_S=\int_{0}^{\infty} p(x)q(x)s(x)x^{\alpha
-1}e^{-x}dx+\lambda \int_{0}^{\infty} p'(x)q'(x)x^{\alpha
}e^{-x}dx\, ,$$ where $s(x)=(1+\delta)x-\alpha\, \delta.$
Therefore, the inner product associated with the measures in
(\ref{ltype11}) is a particular case of (\ref{iplsnd}),i.e, for
$\alpha>0$ and $\delta \in (-1,0]$ the inner product
\begin{align}(p,q)_{S*} &:=\frac{1}{(1+\delta)}(p,q)_S=\int_{0}^{\infty}
p(x)q(x)\left(x-\frac{\alpha\, \delta}{1+\delta}\right)x^{\alpha
-1}e^{-x}dx \nonumber\\
&+\frac{\lambda}{1+\delta} \int_{0}^{\infty}
p'(x)q'(x)x^{\alpha }e^{-x}dx :=\int_{0}^{\infty}
p(x)q(x)(x+a_1)x^{\alpha -1}e^{-x}dx \nonumber \\
&+\lambda_1 \int_{0}^{\infty} p'(x)q'(x)x^{\alpha }e^{-x}dx,
\label{acl-lag}
\end{align} where $a_1, \lambda_1\in \mathbb{R}^+ ,$ corresponds
to the case I.1.

The asymptotic behaviour of these polynomials $Q_n^{(ND)}(x)$ and
their scaled $Q_n^{(ND)}(nx)$ is obtained uniformly on compact
subsets of  $\mathbb{C}\setminus [0,\infty)$ and
$\mathbb{C}\setminus [0,4],$ respectively, i.e.

\begin{theorem}[Th. 2 and 4, \cite{mar-mb}]  \label{asin-la-type-1} \hspace{3cm} \begin{enumerate}
\item[(a)]
 Uniformly on compact subsets of $\mathbb{C}\setminus
[0,\infty) \, ,$

\begin{equation} \label{asin-lag-sob}
\lim_{n\to \infty} \frac{Q_n^{(ND)}(x)}{L_{n}^{(\alpha -1)}(x) }=\frac{2(1+\delta)}{%
\sqrt{\lambda^2+4\lambda(1+\delta)} -\lambda
}=\frac{1}{1-(1+\delta)\ell}\, . \end{equation} \item[(b)] For
$\alpha
>0\, ,$ uniformly on compact subsets of $\mathbb{C}\setminus [0,4] $
\begin{equation}
\lim_{n\to \infty}\frac{Q_n^{(ND)}(nx)}{L_{n}^{(\alpha -1 )}(nx)}=
\frac{\displaystyle{ \varphi \left( \frac{x-2}{2} \right) } }{
\displaystyle{ \varphi \left( \frac{x-2}{2} \right) +\ell(1+\delta
) } }. \label{zar}
\end{equation}
 Furthermore,  for $\alpha >-1\, ,$
\begin{equation}
\lim_{n\to \infty}\frac{Q_n^{(ND)}(nx)}{L_{n}^{(\alpha )}(nx)}=
 \frac{\displaystyle{
\varphi \left( \frac{x-2}{2} \right) +1 } }{ \displaystyle{
\varphi \left( \frac{x-2}{2} \right) +\ell(1+\delta ) } },
\label{zar1}
\end{equation}
where \begin{equation}\varphi (x)=x+\sqrt{x^2-1}\,
\label{zhu}\end{equation} with $\sqrt{x^2-1}>0 $ if $x>1\,$, and
$\ell$ is given by \begin{equation} \label{ell-gen}
\ell=\frac{\lambda + 2(1+\delta )-\sqrt{\lambda^2+4\lambda
(1+\delta) }}{2(1+\delta)^2}\in (0,1).\end{equation}
\end{enumerate}
\end{theorem}

\textbf{Remarks.} \begin{itemize} \item  Notice that
(\ref{asin-lag-sob}) was obtained in \cite[Th.3.5]{mei-per-pin}
for the case I.1 if we consider the inner product (\ref{acl-lag}).
Furthermore, there (\ref{asin-lag-sob}) is also established for
the case I.2 given in (\ref{ltype12}).

\item What are the keys to obtain these results? Indeed, two
properties are essential: one is algebraic and the other one is
analytic.
\begin{enumerate}
\item We have the following relation between the standard
polynomials and Sobolev polynomials (see \cite[Lemma 1]{mar-mb})
\begin{equation}
L_{n+1}^{(\alpha -1)}(x)=Q_{n+1}^{(ND)}(x)+a_n\, Q_n^{(ND)}(x)\,
,\quad n\ge 0\, . \label{rel-alg-lag}
\end{equation}
\item The behaviour of the sequence $a_n$ is related to the ratio
$$\frac{\int_0^{\infty} (L_n^{(\alpha)}(x))^2x^{\alpha}e^{-x}dx}{(Q_n^{(ND)},Q_n^{(ND)})_S}$$
when $n$ tends to infinity. Using a non-linear recurrence relation
for this ratio (see (17) in \cite{mar-mb}) and Poincar\'e's
Theorem the asymptotic result was given in \cite[Prop.2]{mar-mb}.
\end{enumerate}
Then, (\ref{asin-lag-sob}) can be deduced. In order to prove
(\ref{zar}--\ref{zar1}) we need to solve the connection problem
for the Laguerre--Sobolev orthogonal polynomials in terms of the
Laguerre polynomials, i.e., using (\ref{rel-alg-lag}) in a
recursive way we get
$$  Q_n^{(ND)}(x)=\sum_{j=0}^{n}\gamma_{j}^{(n)} L_{n-j}^{(\alpha-1)}(x),
\quad n\ge 0, $$ and $\gamma_{0}^{(n)}=1.$ Since the asymptotic
behaviour of the coefficients $\gamma_{j}^{(n)}$ can be obtained
from the asymptotic behaviour of the sequence $\{a_n\}, $ then
(\ref{zar}--\ref{zar1}) can be deduced using the asymptotics of
scaled Laguerre polynomials $L_n^{(\alpha -1)}(nx),$
$L_n^{(\alpha)}(nx), $ and Lebesgue's dominated convergence
theorem. Although case I.2 was not considered in \cite{mar-mb},
the same techniques used there allow us to deduce (\ref{zar1}) for
case I.2.

\item In \cite[Th.4]{mar-mb}  to get (\ref{zar}--\ref{zar1}) it is
assumed that $\delta \in (-1,0].$ This condition was posed in
order to guarantee a dominant for the sequence \-
$Q_n^{(ND)}(nx)/L_n^{(\alpha -1)}(nx)\, $ on compact subsets of
$\mathbb{C}\setminus [0,4].$ However, this restriction can be
removed by adapting Lemma 3.2 in \cite{alf-mb-rez} to our
circumstances, and then we get a dominant term for
$Q_n^{(ND)}(nx)/L_n^{(\alpha -1)}(nx)$ for any value of $\delta.$

\item From this theorem we can obtain the outer strong asymptotic
of the Laguerre--Sobolev orthogonal polynomials (taking into
account the Perron's Theorem) and scaled ones (from a result by
J.S. Geronimo and W. Van Assche \cite[formula
(4.2)]{waltergeronimo}).

\end{itemize}

\noindent In \cite{mar-mb} the asymptotics of the sequence
$Q_n^{(ND)}$ on the oscillatory region is obtained:

\begin{theorem}[Th. 3, \cite{mar-mb}]
\label{asint-interior-lag}  \item Uniformly on compact subsets of
$(0, \infty)\, ,$
$$\frac{Q_n^{(ND)}(x)}{n^{(\alpha-1)/2}}=e^{x/2}\,
x^{-(\alpha-1)/2} J_{\alpha-1}(2\sqrt{n x})+O(n^{-1/4})\, ,$$
where $ J_{\alpha-1}(x)$ is the Bessel function
$$ J_{\alpha-1}(x)=\sum_{j=0}^{\infty}\frac{(-1)^j (x/2)^{2
j+\alpha-1}}{j! \Gamma (j+\alpha)} \, .$$
\end{theorem}

\noindent \textbf{Remarks.}

\begin{enumerate}

\item To prove the above Theorem the Fej\'{e}r's formula for
Laguerre polynomials (see \cite[Th.8.22.1, p.198]{sz}) is needed
as well as an asymptotic result for the Bessel function
$J_{\alpha-1}(x), $ together with the asymptotics of the sequence
$a_n.$

\item From Theorems \ref{asin-la-type-1} and
\ref{asint-interior-lag} we can deduce that the non--diagonal
Laguerre--Sobolev polynomials $Q_n^{(ND)}$ and, therefore the
Laguerre-Sobolev polynomials of case I.1,  are asymptotically
equivalent to the Laguerre ones of parameter $\alpha-1.$
\end{enumerate}

\subsubsection{Laguerre--Sobolev orthogonal polynomials of case II}

Now we introduce the asymptotic behaviour of the polynomials
$Q_n^{(L2)}(x)=\left((-1)^n/n!\right) x^n+\ldots$ orthogonal with
respect to a Sobolev inner product involving the measures given in
(\ref{ltype2}), i.e,
$$(p,q)=\int_{0}^{\infty} p(x)\, q(x)\, x^{\alpha} e^{-x}dx+\lambda \int_0^{\infty}
p'(x)\, q'(x)\, \frac{x^{\alpha+1} e^{-x}}{x+a}dx+ Mp'(-a)\,
q'(-a),$$ where $a\ge 0, \, \, \alpha
>-1$ and $ M\ge 0.$

The first result appears in \cite{mei-per-pin}:

 \begin{theorem}[Th. 4.11, \cite{mei-per-pin}] \label{lmei} \item
 Relative outer asymptotics.
\begin{equation}
\lim_{n\to \infty}
\frac{Q_n^{(L2)}(x)}{L_n^{(\alpha-1)}(x)}=\left\{
\begin{array}{l}
\displaystyle{\frac{1}{1-\ell}}\left( 1+\frac{a^{1/2}}{(-x)^{1/2}}
\right) \quad
\mathrm{if} \quad M=0, \\
\displaystyle{\frac{1}{1-\ell}}\left( 1-\frac{a^{1/2}}{(-x)^{1/2}}
\right) \quad \mathrm{if} \quad M>0, \label{asi-lag-2}
\end{array}
\right.
\end{equation}
uniformly on compact subsets of $ \mathbb{C}\setminus [0, \infty),
$ where
\begin{equation} \label{ele} \ell=\frac{\lambda+2-\sqrt{\lambda^2+4\lambda}}{2}\in
(0,1).\end{equation}
\end{theorem}

\noindent \textbf{Remark.} Notice that (\ref{ell-gen}) is
(\ref{ele}) when we consider the inner product (\ref{iplsnd}) as
\textit{diagonal}, that is,  $\delta=0.$

In \cite{alf-mb-rez} it was proved that the Sobolev polynomials of
case II have the same outer scaled asymptotic and the same
asymptotics on the oscillatory region as in the case I. It is
convenient to remember the expression of the inner product for
Laguerre--Sobolev of type I in terms of non--diagonal
Laguerre--Sobolev given by (\ref{acl-lag}).

\begin{theorem}[Th. 3.1 and 3.3(b) in \cite{alf-mb-rez}] \label{lasin2}
 \item Theorem \ref{asint-interior-lag} and (\ref{zar1}) in Theorem \ref{asin-la-type-1}
  hold for the polynomials $Q_n^{(L2)}(x).$
\end{theorem}

\noindent \textbf{Remarks.}
\begin{itemize}
\item The key idea of the proofs of Theorems \ref{lmei},
\ref{lasin2} is the same as in the proof of Theorem
\ref{asin-la-type-1} but
 more complicated from a technical point of view. Now, we have an algebraic
relation:
$$L_n^{(\alpha)}(x)-\sigma_{n-1}
L_{n-1}^{(\alpha)}(x)=Q_n^{(L2)}(x)-a_{n-1}^{(L2)}
Q_{n-1}^{(L2)}(x)$$ and the analysis of the asymptotic behaviour
of the sequence $\sigma_n$ is not trivial (see
\cite{mei-per-pin}). Furthermore, as it was shown in
\cite{alf-mb-rez}, it is necessary to study the analytic
\label{tes} properties of the polynomials
$V_n(x):=L_n^{(\alpha)}(x)-\sigma_{n-1} L_{n-1}^{(\alpha)}(x)$ and
one gets $V'_n(x)=-T_{n-1}(x)$ where
$T_{n}(x)=\left((-1)^n/n!\right)\, x^n+\ldots $ are the
polynomials orthogonal with respect to the standard inner product:
$$(p,q)=\int_{0}^{\infty} p(x)\, q(x)\, \frac{x^{\alpha+1} e^{-x}}{x+a}dx+
Mp(-a)q(-a). $$ This family $\{T_n\}$ is orthogonal with respect
to a measure whose absolutely continuous part is a rational
modification of the Laguerre weight function $x^{\alpha+1}e^{-x}$
on $[0, \infty)$ and possibly with a mass point (a Dirac delta) at
$-a\le 0.$ \item \label{pregunta} One important question arises:
\textit{why (\ref{asi-lag-2})  does not yield a constant as in
case I?} We will come back to this question later on.

\end{itemize}

\subsubsection{Zeros and its asymptotics}

In \cite{bru-mei} a complete study of the zeros of coherent pairs
of Jacobi and Laguerre types using Gaussian quadrature formulas
was done. Before stating the results  it is important to point out
that the notation used in  \cite{bru-mei} for  Laguerre--Sobolev
of case I
 is: $$ d\mu_0=(x+a)x^{\alpha} e^{-x}dx,\,
d\mu_1=x^{\alpha+1} e^{-x}dx, \quad a\ge 0, \, \, \alpha
>-1,
$$
and so, the change $\alpha\rightarrow \alpha-1$ must be made in
the results of \cite{bru-mei} for case I to agree with  the
notation (\ref{ltype11}).

\begin{theorem}[Th. 4.1, 4.2, 4.3, 4.4. and 4.5 in \cite{bru-mei}] \label{ceros-lag-1}
\item  Let $Q_n^{(L1)}$ be the La\-guerre-Sobolev orthogonal
polynomials of case I, i.e., orthogonal with respect  to Sobolev
inner products induced by the measures given by  (\ref{ltype11})
or (\ref{ltype12}), then for $n\ge 2,$
\begin{itemize}
\item All the zeros of $Q_n^{(L1)}$ are real, simple, and
positive. \item The zeros of $Q_n^{(L1)}$ interlace with the zeros
of
 $L_n^{(\alpha)},$  i.e., denoting by $l_{n,i}^{(\alpha)}$
and $q_{n,i}^{(L1)}$ the zeros of $L_n^{(\alpha)}$ and
$Q_n^{(L1)},$ respectively, then we get
\begin{equation} \label{inter-lag-1} q_{n,1}^{(L1)}<l_{n,1}^{(\alpha)}<q_{n,2}^{(L1)}<l_{n,2}^{(\alpha)}<\ldots
<q_{n,n}^{(L1)}<l_{n,n}^{(\alpha)}\, .\end{equation} \item For
case I.1, the zeros of $Q_n^{(L1)}$ interlace with the zeros of
 $L_n^{(\alpha-1)},$  i.e.,
\begin{equation} \label{inter-lag-2} l_{n,1}^{(\alpha-1)}<q_{n,1}^{(L1)}<l_{n,2}^{(\alpha-1)}<q_{n,2}^{(L1)}
<\ldots <l_{n,n}^{(\alpha-1)}<q_{n,n}^{(L1)}\, \end{equation} and
for case I.2, the zeros of $Q_n^{(L1)}$ interlace with  the zeros
of $xL_{n-1}^{(1)}$ as in (\ref{inter-lag-2}). \item The zeros of
$Q_{n-1}^{(L1)}$ separate the zeros of $Q_n^{(L1)}.$
\end{itemize}
\end{theorem}

For case II the results are summarized in the following theorem:
\begin{theorem}[Th. 4.1 in \cite{bru-mei}] \label{zerosl2}
\item Let $Q_n^{(L2)}$ be the Laguerre-Sobolev orthogonal
polynomials of case II, then for $n\ge 2,$
\begin{itemize}
\item All the zeros of $Q_n^{(L2)}$ are real, simple, and at most
one of them is in $(-\infty, 0].$ \item The zeros of $Q_n^{(L2)}$
interlace with the zeros of Laguerre $L_n^{(\alpha)}$ in the same
sense as (\ref{inter-lag-1}).
\end{itemize}
\end{theorem}

\noindent \textbf{Remarks.} These results give us bounds for the
zeros of Laguerre--Sobolev polynomials independently of the
parameter $a$ and $\lambda.$ On the other hand, it is important to
note that the other  results appearing in Theorem
\ref{ceros-lag-1} for case I are not true in general for case II.

In order to obtain asymptotic properties of the zeros of
Laguerre--Sobolev polynomials we obtain Mehler--Heine type
formulas (see \cite[p.193]{sz}).

The study of the asymptotics of the zeros of Laguerre--Sobolev
polynomials was done in the papers \cite{mar-mb} for case I.1 and
in \cite{alf-mb-rez} for case II. The results were obtained from
some generalizations of the  Mehler--Heine formula for Laguerre
polynomials.

\begin{theorem}[Th. 5 in \cite{mar-mb} and Th. 3.3 (a) in
\cite{alf-mb-rez}] \label{th-mhl} \hspace{3cm} \begin{enumerate}
\item[(a)] Let $Q_n^{(ND)}$ be the  polynomials orthogonal with
respect to (\ref{iplsnd}) and $\alpha
> -1\, ,$
\begin{equation} \label{mhl-1}
\lim_{n\to \infty} \frac{Q_{n}^{(ND)}(x/n)}{n^{\alpha
-1}}=\frac{1}{1- (1+\delta) \ell }\, x^{-(\alpha-1)/2} J_{\alpha
-1}(2\sqrt{x}) \, .
\end{equation}
holds uniformly for $x$ on compact subsets of $\mathbb{C}\, $ and
$\ell$ is given by (\ref{ell-gen}).

\item[(b)] Let $Q_n^{(L2)} $ be the Laguerre--Sobolev orthogonal
polynomials of case II. If $a>0,$ then
$$\lim_{n\to \infty}\frac{Q_n^{(L2)}(x/n)}{n^{\alpha
-1/2}}=- \frac{d(a)}{1-\ell}\,x^{- \alpha /2}\,J_{\alpha}(2 \sqrt
x).$$ If $a=0$ and $M>0,$ then $$ \lim_{n\to \infty}
\frac{Q_n^{(L2)}(x/n)}{n^{\alpha -1}}=\frac{1}{1-\ell}\,s(x).$$
Finally, if $a=M=0,$ then (\ref{mhl-1}) holds with $\delta=0.$

\noindent For
  $\ell$  given by (\ref{ele}),
$$d(a)=\begin{cases} -\sqrt{a}& \textrm{if} \quad
M=0, \\ \sqrt{a}& \textrm{if} \quad  M
>0,
\end{cases} $$
 and
\begin{equation} \label{mh-vn-esp}
s(x)=x^{-(\alpha-1)/2}J_{\alpha-1}(2
\sqrt{x})-(\alpha+1)\,x^{-\alpha/2}J_{\alpha}(2 \sqrt{x}).
\end{equation}
All the limits hold uniformly on compact subsets of $\mathbb{C}. $
\end{enumerate}

\end{theorem}

\noindent \textbf{Remarks.} As we have explained in  the third
remark of Theorem \ref{asin-la-type-1}, the restriction on
$\delta$ that appears in \cite{mar-mb} can be removed. On the
other hand, (\ref{mhl-1}) holds for the Laguerre--Sobolev
polynomials of case I:  it is clear for case I.1 as we already
know and it follows in a straightforward way for case I.2.

Now, we can obtain a limit relation between the zeros of
Laguerre--Sobolev polynomials and the zeros of Bessel functions
using Hurwitz's theorem in the results of this theorem. Thus we
get

\begin{corollary}[Prop. 4 in \cite{mar-mb} and Cor. 3.4
in \cite{alf-mb-rez}] \hspace{3cm}
\begin{enumerate} \item[(a)]  Let  $q_{n,i}^{(ND)}$ be the
zeros of $Q_n^{(ND)}$. Then \begin{equation} \label{ceros-lag1}
\lim_{n\to \infty} n\,q_{n,i}^{(ND)}=\frac{j_{\alpha-1,i}^2}{4} \,
,\end{equation} where three cases appear: \begin{itemize} \item If
$-1<\alpha<0,$  $j_{\alpha-1,1}$ is any of the two purely
imaginary zeros of $J_{\alpha-1}(x)$ and, for $i\ge 2,$
$j_{\alpha-1,i}$ is the $(i-1)$th positive real zero of
$J_{\alpha-1}(x).$ \item If $\alpha=0,$
$j_{\alpha-1,1}=j_{-1,1}=0$  and, for $i\ge 2,$ $j_{-1,i}$ is the
$(i-1)$th positive real zero of $J_{-1}(x).$ \item If $\alpha>0,$
$j_{\alpha-1,i}$ is the $i$th positive real zero of
$J_{\alpha-1}(x).$
\end{itemize}
\item[(b)] Let  $q_{n,i}^{(L2)}$ be the zeros of $Q_n^{(L2)}$.
 \begin{itemize} \item If $a>0,$ then we have
$$\lim_{n\to \infty} n\,q_{n,i}^{(L2)}=\frac{j_{\alpha,i}^2}{4},
$$ where $j_{\alpha,i}$ the ith positive zero
of the Bessel function $J_\alpha(x).$

\item If $a=0$ and $M>0,$ then we have
 $$
 \lim_{n\to \infty} n\,q_{n,i}^{(L2)}=s_{\alpha,i},$$
where $s_{\alpha,i}$ denotes the $i$th real zero of function
$s(x)$ defined in (\ref{mh-vn-esp}).

\item If $a=M=0,$ then (\ref{ceros-lag1}) holds.

\end{itemize}

\end{enumerate}

\end{corollary}

\noindent \textbf{Remarks.}\begin{itemize} \item  Notice that  for
the zeros of Laguerre--Sobolev polynomials of case I,
(\ref{ceros-lag1}) holds. Indeed,  cases I.1 and I.2 correspond to
the third and second case, respectively. \item The
Laguerre--Sobolev polynomials of case II have exactly one negative
zero for $n$ sufficiently large if and only if either $ \alpha>-1,
\, a=0, \, \mathrm{and} \, M>0$ or $ -1<\alpha<0 \, \,
\mathrm{and} \,\,  a=M=0.$

\item For a fixed $n,$ it would be very interesting to analyze the
behaviour of $q_{n,j}^{(Lk)}\, $ $\, k=1,2\, $ in terms of the
parameter $\lambda .$

\item The function $s(x)$ given by (\ref{mh-vn-esp}) has only one
negative zero and their positive zeros separate  those of
$J_{\alpha+1}(2\sqrt{x}).$

 \item From
this theorem we can get information about the critical points. In
fact, the critical points of these families of polynomials  for
$n$ sufficiently large lie on $[0,+\infty).$

\end{itemize}

\subsubsection{A generalization of Laguerre polynomials: the Laguerre--Sobolev type}

In \cite{koe-mei}  the inner product
\begin{equation} \label{prmn}
(p,q)=\frac{1}{\Gamma(\alpha+1)}\int_{0}^{\infty} p(x)
q(x)\,x^{\alpha}\, e^{-x}dx+M\, p(0)q(0)+ N\, p'(0) q'(0)\, ,
\end{equation}
where $M, N \ge 0$ and $\alpha>-1$ was considered. The polynomials
$\big(L_n^{(\alpha,M,N)}(x)\big)_n$ orthogonal with respect to
(\ref{prmn}) constitute a natural extension of the so--called
Koornwinder's generalized Laguerre polynomials which are
orthogonal with respect to (\ref{prmn}) where $N=0$. Algebraic and
differential properties of the generalized Laguerre polynomials
$\big(L_n^{(\alpha,M,N)}(x)\big)_n$ have been obtained in the 90's
of the past century. However, the asymptotic properties of
$\big(L_n^{(\alpha,M,N)}(x)\big)_n$ were not known until its
consideration in \cite{ran-mb}.

The asymptotic of these polynomials is obtained  in
$\mathbb{C}\setminus [0, \infty)$ as well as in the oscillatory
region $(0,\infty).$ In \cite{ran-mb}, following the notation of
\cite{koe-mei}, $L_n^{(\alpha)}(x)$ and $L_n^{(\alpha,M,N)}(x)$
have as leading coefficients, respectively $(-1)^n/n!$  and
$$\frac{(-1)^n}{n!}\, \left(B_0(n)-nB_1(n)+n(n-1) B_2(n)\right),
$$  where
\begin{align*} B_0(n) &=
1-\frac{N}{\alpha\!+\!1}{n\!+\!\alpha\!+\!1\choose n\!-\!2}
,\,\,\, B_1(n) =-\frac{M}{\alpha\!+\!1} { n\!+\!\alpha\choose
n}-\frac{(\alpha\!+\!2) N}
{(\alpha\!+\!1)(\alpha\!+\!3)} {n\!+\!\alpha\choose n\!-\!2}, 
\\[4mm]
B_2(n)&=\frac{N}{(\alpha\!+\!1)(\alpha\!+\!2)(\alpha\!+\!3)}
{n\!+\!\alpha\choose n\!-\!1}+
\frac{MN}{(\alpha\!+\!1)^2(\alpha\!+\!2)(\alpha\!+\!3)}
{n\!+\!\alpha\choose n}{n\!+\!\alpha\!+\!1\choose n\!-\!1}\, .
\end{align*}

But if we take $(-1)^n/n!$ as leading coefficient in both families
of polynomials  then the computations can be simplified. For
instance,
 in the following theorem  we summarize some of the results obtained
in \cite{ran-mb}:

\begin{theorem}[Th. 1 and 2(a) in \cite{ran-mb}] \hspace{3cm} \begin{enumerate}
\item[(a)] Relative asymptotics. $$ \lim_{n\to
\infty}\frac{L_n^{(\alpha, M,N)}(x)}{L_{n}^{(\alpha)}(x)}=1, \quad
\mathrm{with} \quad M,N\geq 0,\, \, \alpha>-1,
$$ uniformly on compact subsets
of $\mathbb{C}\setminus [0, +\infty).$ \item[(b)] Strong
asymptotics on compact subsets of $(0, +\infty). $

 Denote by
$$g_{n,i}(x)= x^{-\alpha/2}J_{\alpha+2i}\left(2\sqrt{n x}\right),
$$ then
\begin{align*}
&\frac{L_n^{(\alpha, M,N)}(x)}{n^{\alpha/2}}= \\ &\left\{%
\begin{array}{ll}
    c_1(n)e^{x/2}\, g_{n-1,1}(x)+O\left(n^{-\min\{\alpha+5/4,3/4\}}\right), & \hbox{$M>0$, N=0;}
    \\ \\
\begin{array}{l}
    c_2(n)e^{x/2}\left( a_2(n)g_{n-2,2}(x)-(\alpha+2)a_1(n) g_{n-1,1}(x)\right.
    \\ \left. -a_0(n)g_{n,0}(x)\right)+O\left(n^{-3/4}\right),
    \end{array}
     & \hbox{M=0, $N>0$;} \\ \\
    c_3(n)e^{x/2}\, g_{n-2,2}(x)+O\left(n^{-\min\{\alpha+5/4,3/4\}}\right), & \hbox{M, N$>$0,} \\
\end{array}%
\right.     \end{align*} where $\lim_{n\to \infty}c_i(n)=1,$
$i=1,3,$ $\lim_{n\to \infty}c_2(n)=1/(\alpha+2), $ and $\lim_{n\to
\infty}a_i(n)=1, $ $i=0,1,2.$

\item[(c)] Relative outer asymptotics of scaled polynomials.

  $$\lim_{n\to
\infty}\frac{L_n^{(\alpha, M,N)}(nx)}{L_{n}^{(\alpha)}(nx)}=1,
\quad \mathrm{with} \quad M,N\geq 0,\, \, \alpha>-1,
$$ uniformly on compact subsets of $\mathbb{C}\setminus [0,4]\, .$

\end{enumerate}

\end{theorem}

\noindent \textbf{Remark.} It is important to note that we have
quite simplified the results appearing in \cite{ran-mb} as well as
we have corrected some misprints, in particular one of them in the
scaled asymptotics for the case $M=0, $ $N>0.$ For all of them,
now it is more easy to understand the meaning of the results.

\bigskip

Concerning their zeros, in \cite{koe-mei} it was proved that the
zeros of polynomials $L_n^{(\alpha, M,N)}(x)$ are real and at
least $n-1$ of them lie in $(0,\infty).$ Furthermore, when $N>0$
and $n$ is large enough there exists exactly one zero in
$(-\infty, 0].$ We denote by $q_{n,i}$ the zeros of $L_n^{(\alpha,
M,N)}(x).$ Now, let $N>0$ and $n$ sufficiently large such that
$L_n^{(\alpha, M,N)}(x)$ has the zero $q_{n,1}$ in $(-\infty,0].$
Then, for $M>0,$ $$ -\frac{1}{2}\, \left( \frac{N}{M}
\right)^{1/2}\leq q_{n,1}\leq 0.
$$
It would be very interesting to analyze the behaviour of $q_{n,1}$
in terms of $N,$ i.e., $q_{n,1}=O(N^{-\mu})$ with $\mu \in
\mathbb{R}^{+},$ as well as, the behaviour of $q_{n,j}$ in terms
of $N,$ for a fixed $M.$ Notice that for $N=0,$ a standard case,
$q_{n,1}=O(M^{-1}),$ and $q_{n,j}$ is a decreasing function in
terms of $M.$

 Later on in \cite{ran-mb} asymptotic properties
of the zeros were looked for. So Mehler--Heine type formulas were
obtained and therefore limit relations between the zeros of
generalized Laguerre polynomials and the zeros of Bessel functions
were also deduced. We summarize the results as:
\begin{proposition}[Th. 2 (b) and 3 in \cite{ran-mb}]
 \item If we denote
\begin{equation} \label{fun-mh-g} g_{i}(x)=
x^{-\alpha/2}J_{\alpha+2i}\left(2\sqrt{x}\right),
\end{equation} then
$$
 \lim_{n\to
\infty}\frac{L_n^{(\alpha, M,N)}(x/n)}{n^{\alpha}}=\left\{%
\begin{array}{ll}
    -g_1(x), & \hbox{$M>0$, N=0;}
    \\ \\
\frac{1}{\alpha+2}\left( g_2(x)-(\alpha+2)g_1(x)-g_0(x)\right)
     & \hbox{M=0, $N>0$;} \\ \\
    g_2(x), & \hbox{M, N$>$0.} \\
\end{array}%
\right. $$ uniformly on compact subsets of $\mathbb{C}.$
Therefore, denoting by $j_{\alpha,i}$ the $i$-th positive zero of
the Bessel function $J_\alpha(x)$ and by $q_{n,i}$  the zeros in
increasing order of the polynomial $L_n^{(\alpha, M,N)}(x)$
 with
$\alpha>-1,$
\begin{enumerate}
\item[(a)] If $M>0$ and $N=0,$ then $$ \lim_{n\to \infty}
n\,q_{n,1}=0 \quad \mathrm{and}\quad \lim_{n\to \infty}
n\,q_{n,i}=\frac{j_{\alpha+2,i-1}^2}{4}\, , \, i\ge 2.
$$

\item[(b)] If $M=0$ and $N>0,$ then
$$
 \lim_{n\to \infty} n\,q_{n,i}=z_{i},
$$
where $z_{i}$ denotes the $i$-th real zero of the function
$g(x)=g_2(x)-(\alpha+2)g_1(x)-g_0(x)$ with $g_i(x)$ given by
(\ref{fun-mh-g}). Moreover, $g(x)$ has only one negative real
zero.

\item[(c)] If $M,N>0,$ then
$$
\lim_{n\to \infty} n\,q_{n,i}=0, \, i=1,2 \quad \mathrm{and} \quad
\lim_{n\to \infty} n\,q_{n,i}=\frac{j_{\alpha+4,i-2}^2}{4}\, , \,
\, i\ge 3. $$
\end{enumerate}
\end{proposition}

\noindent \textbf{Remarks.} Notice that for $M>0,$ $N=0$ the first
scaled zero, $nq_{n,1},$ of the orthogonal polynomials
(Koornwinder polynomials) $L_n^{(\alpha, M,0)}(x)$ tends to $0$
for $n\to \infty$. On the other hand, in the case $M>0,$ $N>0$ the
two first scaled zeros $nq_{n,i},$ $i=1,2$, with $q_{n,1}<0$,
also tend to $0$ when $n\to \infty.$ In both cases the results
agree with those of \cite{alr} in a more general framework.
However,  in the case  $M=0,$ $N>0$  for $n$ large enough, the
first scaled zero $nq_{n,1}$ is always negative and does not tend
to zero when $n\to\infty$, i.e., $\lim_{n\to
\infty}nq_{n,1}=z_{1}$ is a negative real number.

In \cite{ran-mb} a conjecture was stated for the polynomials
$L_n^{(\alpha,M_0,\ldots, M_s)}(x)$  orthogonal with respect to
the inner product
$$
(p,q)=\frac{1}{\Gamma(\alpha+1)}\int_{0}^{\infty} p(x)
q(x)\,x^{\alpha}\, e^{-x}dx+\sum_{i=0}^{s}M_i\,
p^{(i)}(0)q^{(i)}(0)\, , \quad \alpha>-1.
$$
Now, after the simplifications made in the results of
\cite{ran-mb},  we can formulate it  as follows:

\noindent \textit{\underline{Conjecture.}} If $M_i>0,$
$i=0,\ldots, s\, $ with $s\ge 2$, then $$
\lim_{n\to\infty}\frac{L_n^{(\alpha, M_0,\ldots,
M_s)}(x/n)}{n^{\alpha}}=(-1)^{s+1}\, g_{s+1}(x)\, ,
$$
uniformly on compact subsets of $\mathbb{C}$ where $g_{s+1}(x)$ is
given by (\ref{fun-mh-g}).

From the above proposition we can not formulate this conjecture
when some of the constants $M_i,\,$ $i=0,\ldots, s-1,$ vanish with
$M_s>0.$ What occurs in this case? This is an interesting unsolved
problem.

\subsection{Hermite weights}
Next we consider the inner product
\begin{equation} \label{iphermite}(p,q)_S=\int_{-\infty}^{\infty} p(x)q(x)d\mu_0+\lambda
\int_{-\infty}^{\infty} p'(x)q'(x)d\mu_1\end{equation} where
$(\mu_0,\mu_1)$ is a symmetrically coherent pair of measures (see
\cite{meijer}) corresponding to Hermite weight, that is, the
measures are given by (\ref{htype1}) (case I) or (\ref{htype2})
(case II) according to the classification given by Meijer, i.e,
$$\begin{array}{lll}d\mu_0=(x^2+a^2) e^{-x^2}dx,& d\mu_1=
e^{-x^2}dx,& a\in \mathbb{R}, \quad
\mathrm{or} \\
 d\mu_0=e^{-x^2}dx,& d\mu_1=e^{-x^2}/(x^2+a^2)dx,& a\in \mathbb{R}\setminus \{0\}.
 \end{array} $$

In this section we use the notation: $H_n(x)=2^n x^n+\ldots$ for
the classical Hermite polynomials and $Q_n^{(H1)}(x)=2^n
x^n+\ldots\, ,$ $\, Q_n^{(H2)}(x)=2^n x^n+\ldots$ for the Sobolev
polynomials orthogonal with respect to (\ref{iphermite}) where the
measures are given by (\ref{htype1}) and (\ref{htype2}),
respectively.

In the framework of symmetrically coherent pairs, we have the
algebraic relation between the classical Hermite polynomials and
the Sobolev orthogonal polynomials (see, for example,
\cite[Th.4(c)]{iserles}) $$ H_{n+2}(x)+\sigma_{n+1}^{(Hi)}\,
H_n(x)=Q_{n+2}^{(Hi)}(x)+a_n^{(Hi)}Q_n^{(Hi)}(x), \quad n\ge i-1,
\quad i=1,2,$$ where $\sigma_{n}^{(H1)}=0$ for all $n\in
\mathbb{N}.$

As in the Laguerre case the keys to obtain the asymptotic
behaviour of the Sobolev orthogonal polynomials are related to the
limits of the sequences $\sigma_{n}^{(H2)}$ and $a_n^{(Hi)},$
$i=1,2,$ when $n\to \infty.$ It is important to note that the
sequences $a_n^{(H1)}$ and $a_n^{(H2)}$ involve the ratios $$
\frac{\int_{-\infty}^{\infty}
H_n^2(x)e^{-x^2}dx}{(Q_n^{(H1)},Q_n^{(H1)})_S}, $$
$$\frac{\int_{-\infty}^{\infty}
 W_n^2 (x)\,
\frac{e^{-x^2}}{x^2+a^2}dx}{(Q_n^{(H2)},Q_n^{(H2)})_S},$$
respectively where $W_n(x)=2^nx^n+\ldots$ are the polynomials
orthogonal with respect to the weight function
\begin{equation} \label{ip-1-h} \frac{e^{-x^2}}{x^2+a^2}, \quad \quad  a\neq 0.
\end{equation}
  From a technical point of view the case I
is  more simple   than case II. In fact, in \cite[Lemma 2.2, 2.4
and 2.6]{alf-mb-per-pin-rez} we get:
\begin{equation} \label{lim-par-h-1}
\lim_{n\to \infty}\frac{\sigma_{n}^{(H2)}}{2n}=1, \quad \lim_{n\to
\infty}\frac{a_{n}^{(Hi)}}{2n}=\frac{1}{\varphi(2\lambda+1)},\quad
i=1,2, \end{equation} where $\varphi$ is given by (\ref{zhu}).
This is enough to deduce the relative asymptotics between Sobolev
orthogonal polynomials and Hermite polynomials in case I. However,
for case II  more information about the asymptotics of
$\sigma_{n}^{(H2)}/(2n)$ is needed to get an useful result. So, in
\cite[Prop. 2]{mb} it was deduced
\begin{equation} \label{lim-par-h-2}\lim_{n\to \infty}
\sqrt{[n/2]}\left(\frac{\sigma_{n}^{(H2)}}{2n}-1\right)=-|a|\,
.\end{equation} From it  we get:
\begin{theorem}[Th. 2.3 in \cite{alf-mb-per-pin-rez} and Th. 2 in
\cite{mb}] \item Let
\begin{equation} \label{theta}
\theta(\lambda)=\frac{\varphi(2\lambda+1)}{\varphi(2\lambda+1)-1}\,
.
\end{equation}
\begin{itemize}
\item[(a)] Outer relative asymptotics.
\begin{equation*}
\lim_n\frac{Q_n^{(H1)}(x)}{H_n(x)}=\theta(\lambda)\, ,
\end{equation*}
\noindent uniformly on compact subsets of $\mathbb{C} \setminus
\mathbb{R}$. \item[(b)] Outer relative asymptotics.
\begin{equation} \label{asi-her-2}\lim_{n\to
\infty}\sqrt{\left[\frac{n}{2}\right]}\,
\frac{Q_n^{(H2)}(x)}{H_n(x)}=\left\{ \begin{array}{cll}
\displaystyle{\theta(\lambda)\,  \left(-i\,
x+|a|\right)} & \mathrm{if} & x\in \mathbf{C}_+, \\
\displaystyle{\theta(\lambda)\,\left(i\, x+|a|\right)} &
\mathrm{if} & x\in \mathbf{C}_-, \end{array} \right.\end{equation}
uniformly on compact subsets of half planes
$\mathbf{C}_+:=\{x=\alpha+i\, \beta \in \mathbb{C}\, : \beta>0\}$
and $\mathbf{C}_-:=\{x=\alpha+i\, \beta \in \mathbb{C}\, :
\beta<0\}$, respectively.
\end{itemize}
Again, $\varphi$  is given by (\ref{zhu}).
\end{theorem}

\noindent \textbf{Remark.} \label{pregunta-2} As in the Laguerre
case, in the relative asymptotic we get as limit a constant  for
case I and a function depending on $x$ for case II. Thus,
\textit{why  these different situations appear?}. Later on, we
will come back on this question.

On the other hand, the asymptotic of scaled Hermite polynomials
was obtained in \cite[Th. 3.1]{alf-mb-per-pin-rez}:

\begin{theorem}[Th. 3.1 in \cite{alf-mb-per-pin-rez}] \hspace{3cm}
\begin{itemize} \item[(a)] Outer relative asymptotics for scaled
polynomials.
$$ \lim_n\frac{Q_n^{(H1)}(\sqrt{n}x)}{H_n(\sqrt{n}x)}=\frac{\varphi(2\lambda+1)
\varphi^2\left(\frac{x}{\sqrt
2}\right)}{\varphi(2\lambda+1)\varphi^2\left(\frac{x}{\sqrt
2}\right)+1}
$$
\noindent  uniformly on compact subsets of $\mathbb{C} \setminus
[-\sqrt 2,\sqrt 2]$.

\item[(b)] Outer relative asymptotics for scaled polynomials.
$$ \lim_n\frac{Q_n^{(H2)}(\sqrt{n}x)}{H_n(\sqrt{n}x)}=\frac{\left(\varphi^2\left(\frac{x}{\sqrt
2}\right)+1\right)\varphi(2\lambda+1)}{\varphi(2\lambda+1)\varphi^2\left(\frac{x}{\sqrt
2}\right)+1}
$$
\noindent  uniformly on compact subsets of $\mathbb{C} \setminus
[-\sqrt 2,\sqrt 2]$.
\end{itemize}
Again, $\varphi$  is given by (\ref{zhu}).
\end{theorem}

It is important to note that for the scaled asymptotics of case II
it is enough to know (\ref{lim-par-h-1}). However, for the
non-scaled asymptotics obtained later in \cite{mb} some extra
information about the sequence $\sigma_{n}^{(H2)}$ given by
(\ref{lim-par-h-2}) is needed. This follows from the fact that
\begin{align*}
\lim_{n\to \infty}\frac{H_{n+2}(x)+\sigma_{n+1}^{(H2)}\,
H_n(x)}{H_{n+2}(x)}&=0, \\
\lim_{n\to \infty}\frac{H_{n+2}(\sqrt{n}\,x)+\sigma_{n+1}^{(H2)}\,
H_n(\sqrt{n}\,x)}{H_{n+2}(\sqrt{n}\,x)}&\ne 0,
\end{align*}
uniformly on compact subsets of $\mathbb{C} \setminus \mathbb{R}$
and of $\mathbb{C} \setminus [-\sqrt 2,\sqrt 2],$ respectively.

Next we give the inner asymptotics of the Sobolev polynomials
$S_n^{(Hi)}(x)$ $\, i=1,2$.
\begin{theorem}[Th. 2.3 in \cite{bru-gro-mar-mei-mb}] \label{th-inner}
 Denote
$\lambda_n^*=\displaystyle{\frac{\Gamma(n/2+1)}{\Gamma (n)}}\, .$
\begin{itemize} \item[(a)] For
$n\rightarrow \infty$, the polynomial $Q_n^{(H1)}$ satisfies
\[
\lambda_n^* e^{-\frac{1}{2}x^2} Q_n^{(H1)}(x) = \theta(\lambda)\,
F_n(x) +
 \small{o}(1), \quad x \in \mathbb{R},
\]
where $F_n(x)= \lambda_n^*\, e^{-\frac{1}{2}x^2} H_n(x)$.

\item[(b)] For $n \rightarrow \infty$, the polynomial $Q_n^{(H2)}$
satisfies
\[
 \frac{\lambda_{n-1}^*}{2}\, e^{-\frac{1}{2}x^2} Q_n^{(H2)}(x) =
\theta(\lambda)\, G_n(x) + \small{o}(1), \quad x \in \mathbb{R},
\]
where $G_n(x)=\Big\{\lambda_{n-1}^* xH_{n-1}(x)-|a|\lambda_{n-2}^*
H_{n-2}(x)\Big\} e^{-\frac{1}{2}x^2}$.
\end{itemize}
In both cases the convergence is uniform on compact intervals of
the real line and $\theta(\lambda)$ is given by (\ref{theta}).
\end{theorem}

\noindent \textbf{Remarks.} \begin{itemize} \item Notice that in
\cite[p.200]{sz} the asymptotic behaviour of $F_n$ is given, that
is, for $n\to \infty$ we have
$$
F_n(x) = \cos{\big(\sqrt{2n+1}\, x- \frac{n\pi}{2}\big)} +
\mathcal O \Big(\frac{1}{\sqrt{n}}\Big),
$$
uniformly on compact intervals of the real line. Furthermore, it
is to observe that $G_n(x)=xF_{n-1}(x)-|a|F_{n-2}(x)$ and so, for
$n\to \infty,$ we get
$$G_n(x)=|a| \cos\Big(\sqrt{2n-3}\,x - \frac{n\pi}{2} \Big) -x
\sin\Big( \sqrt{2n-1}\,x - \frac{n\pi}{2} \Big)+
\mathcal{O}\left(\frac{1}{\sqrt{n}} \right),$$ uniformly on
compact intervals of the real line.

\item In \cite[Prop. 2]{cas-mb} some bounds depending of the
degree of the polynomials were obtained  for $|Q_n^{(Hi)}(x)|\, ,$
$i=1,2,$ with $x\in \mathbb{R}.$

\item We also want to emphasize the fact that analytical
properties of Hermite--Sobolev polynomials of case II involve the
study of these properties for the polynomials $W_n$ which are
orthogonal with respect to a rational modification of Hermite
weight function. So, the study of the zeros of the family of
orthogonal polynomials $W_n$ is associated with Gaussian
quadrature formulas for the measure $e^{-x^2}/(x^2+a^2)dx. $

\end{itemize}

\subsubsection{Zeros and its asymptotics}

The zeros of Hermite--Sobolev polynomials were studied in
\cite{bru-gro-mei}. Applying Gaussian quadrature formulas and
other technical tools  the location of the zeros of
Hermite--Sobolev polynomials for both cases I and II is deduced.
We summarize the results as:

\begin{theorem}[Th. 4.7,\, 4.8,\, 4.10,\, 5.7,\, 5.8,\, 5.10 in
\cite{bru-gro-mei}] \label{ceros-hs} \item  Let $h_{n,1}<\ldots
<h_{n,m}$ denote the positive zeros of Hermite polynomials $H_n\,
. $  Then
\begin{itemize}
\item The zeros of $Q_n^{(H1)}$ and $Q_{2n+1}^{(H2)}$ are real and
simple. $Q_{2n}^{(H2)}$ has at least $2n-2$ real and simple zeros.

\item Let $p_{n,1}<\ldots <p_{n,m}$ denote the positive zeros of
polynomials  orthogonal with respect to the weight function
$(x^2+a^2)\,e^{-x^2}$ and let $q_{n,1}^{(H1)}<\ldots
<q_{n,m}^{(H1)}$ denote the positive zeros of  the polynomials
$Q_n^{(H1)}.$ Then, for $n\ge 3,$
\begin{align} &h_{n,1}<q_{n,1}^{(H1)}<\ldots
<h_{n,m}<q_{n,m}^{(H1)}\, , \nonumber
\\
&q_{n,1}^{(H1)}<p_{n,1}<\ldots <q_{n,m}^{(H1)}<p_{n,m}\, .
\label{entre-h1-2}
\end{align}
\item For $n\ge 3,$ the zeros of $Q_{n-1}^{(H1)}$ separate the
zeros of $Q_n^{(H1)}\, .$  \item Let $q_{n,2}^{(H2)}<\ldots
<q_{n,m}^{(H2)}$ denote the $m-1$ largest positive zeros of
polynomials $Q_n^{(H2)}.$ Then,
\begin{align} &h_{n,1}<q_{n,2}^{(H2)}<h_{n,2}<\ldots
<h_{n,m-1}<q_{n,m}^{(H2)}<h_{n,m}\, , \quad n\ge 3,
\label{entre-h2} \\
&h_{n-2,1}<q_{n,2}^{(H2)}<\ldots <h_{n-2,m-1}<q_{n,m}^{(H2)},
\quad n\ge 4. \nonumber
\end{align}
Furthermore, if $Q_n^{(H2)}$ has $n$ real and simple zeros and
denoting by $q_{n,1}^{(H2)}$ the smallest positive zero, then for
$n\ge 3$ we have  $\, q_{n,1}^{(H2)}<h_{n,1}\, .$

\end{itemize}
\end{theorem}

\noindent \textbf{Remarks.}\begin{itemize} \item  In \cite[Th.
5.15]{bru-gro-mei} it is also shown that the zeros of the
polynomials $W_{2n}$ orthogonal with respect to  the weight
(\ref{ip-1-h}) interlace with those of $Q_{2n+1}^{(H2)}$ and the
zeros of $W_{2n+1}$ interlace with those real zeros of
$Q_{2n+2}^{(H2)}.$

\item We can see that this theorem is very similar to Theorems
\ref{ceros-lag-1} and \ref{zerosl2} for the Laguerre case. In
fact, the possibility of one negative zero for Laguerre--Sobolev
polynomials of case II is, in some sense, equivalent to the
possibility of two conjugate complex zeros for even
Hermite--Sobolev
 polynomials of case II.

\item  $Q_{2n}^{(H2)}$ can have complex zeros. For example if we
take $a^2=1/25$ and $\lambda=1, $ then the polynomials
$Q_{2n}^{(H2)}$  have two complex zeros for $n=3,\ldots, 6, $ as
the following table shows
$$
\begin{array}{c|c}
  Q_{2n} & \textrm{Complex   zeros} \\ \hline
  Q_6 & \pm    0.1399747870\, i \\
  Q_{8} & \pm  0.1158504270\, i \\
  Q_{10} & \pm 0.0769761729\, i \\
  Q_{12} & \pm 0.0308661279\, i\\
\end{array}
$$
In \cite{bru-gro-mei} a study of the zeros of $Q_{2n}^{(H2)}$ has
been made when $\lambda\rightarrow \infty.$ In this paper the
authors state that the limit polynomials:
$R_{2n}^{(H2,\infty)}(x):=\lim_{\lambda \to
\infty}Q_{2n}^{(H2)}(x)$ for $n\ge 2,$ have complex zeros when
$|a|\, \sqrt{n-1}<0.21$ and have no complex zeros when $|a|\,
\sqrt{n-1}>2.5$. An open question about the complex zeros of these
polynomials is posed in \cite{bru-gro-mei}: what happens if
$0.21<|a|\, \sqrt{n-1}<2.5\, ?$

\item For a fixed $n,$ and taking into account (\ref{entre-h1-2},
\ref{entre-h2}) we conjecture that $q_{n,k}^{(Hj)}\, $ is an
increasing function of $\lambda$ for each $k.$
\end{itemize}

The study of the asymptotics of the zeros of Hermite--Sobolev
polynomials has been done in the papers \cite{cas-mb} for case I
and in \cite{mb} for case II. As in the Laguerre case, the results have been  obtained from some
generalizations of the  Mehler--Heine type formulas for Hermite
polynomials.

\begin{theorem}[Th. 1 in \cite{cas-mb} and Th. 1 in \cite{mb}]
\item

\begin{enumerate}
\item[(a)] Let  $Q_{n}^{(H1)}$ be the Hermite--Sobolev orthogonal
polynomials of case I. Then,
\begin{align*}
\lim_{n\to \infty}\frac{(-1)^n\sqrt{n}\,Q_{2n}^{(H1)}(x/(2\sqrt{n}))}{
2^{2n}\, n! }&=\theta (\lambda) \,\frac{\cos(x)}{\sqrt{\pi}}\, , \\
 \lim_{n\to \infty}\frac{(-1)^n Q_{2n+1}^{(H1)}(x/(2\sqrt{n}))}
{2^{2n+1}\, n!}&=\theta (\lambda)\,\frac{\sin(x)}{\sqrt{\pi}}\, ,
\end{align*}
hold both  uniformly on compact subsets of $\mathbb{C}\, .$

\item[(b)] Let  $Q_{n}^{(H2)}$ be the Hermite--Sobolev orthogonal
polynomials of case II. Then,
\begin{align*}
\lim_{n\to \infty}\frac{(-1)^n\,Q_{2n}^{(H2)}(x/(2\sqrt{n}))}{
2^{2n}\, (n-1)! }&=\theta (\lambda)
\,\frac{|a|\cos(x)}{\sqrt{\pi}}\, , \\
\quad \lim_{n\to \infty}\frac{(-1)^n
\sqrt{n}\,Q_{2n+1}^{(H2)}(x/(2\sqrt{n}))} {2^{2n+1}\, n!}&=\theta
(\lambda)\,\frac{|a|\sin(x)}{\sqrt{\pi}}\, ,
\end{align*}
hold both uniformly on compact subsets of $\mathbb{C}\, .$
\end{enumerate}
Again, $\varphi$  and $\theta (\lambda)$ are given by (\ref{zhu})
and (\ref{theta}), respectively.
\end{theorem}

Therefore, using Hurwitz's theorem in the results of this theorem
we  obtain the same limit relation between the zeros of
Hermite--Sobolev polynomials for case I and II and the zeros of
the elementary trigonometric functions $\cos(x)$ and $\sin(x).$
That is,

\begin{corollary}[Cor. 1 in \cite{cas-mb} and Cor. 2
in \cite{mb}] \item

 With the notation of Theorem \ref{ceros-hs}, we
have

$$ \lim_{n\to \infty}
2\,\sqrt{n}\,q_{2n,j}^{(Hk)} = (2j-1)\, \frac{\pi}{2}\, ,\quad
\lim_{n\to \infty} 2\, \sqrt{n}\,q_{2n+1,j}^{(Hk)} = j\, \pi\, ,
\quad j\in \mathbb{N}, \quad k=1,2\, .$$

\end{corollary}

\noindent \textbf{Remark.} In \cite{mb} the Mehler--Heine type
formulas  as well as the zero asymptotics for the standard
polynomials  $W_n$ are deduced. In fact, the asymptotic behaviour
of the smallest zeros of the Hermite polynomials $H_n$, the
polynomials $W_n$ and the Hermite--Sobolev orthogonal polynomials
(case I and II) is the same.

  \bigskip

Now is the moment to analyze the question posed in pages
\pageref{pregunta} and \pageref{pregunta-2}: \textit{why the
asymptotic behaviours  for the cases II (Laguerre and Hermite) do
not yield a constant as in the cases I?} The answer is given
taking into account the inner product:
$$(f,g)=\int f(x)g(x)d\mu_0(x)+\lambda\int f'(x)g'(x)d\mu_1(x),
\quad \lambda >0,$$ where $\mu_0, \mu_1$ are positive Borel
measures supported on the real axis. We can see that the
derivative in the second integral introduce a multiplicative
factor $n^2$ when we apply this inner product to monic polynomials
of degree $n.$ Therefore, in ``some sense" the  measure $\mu_1$
plays the most important role and so in terms of relative
asymptotics we should compare the Sobolev polynomials $Q_n$ with
those associated with the second measure $\mu_1.$ This  was
already observed in \cite{alf-and-rez} and it was a motivation to
introduce the Sobolev balanced extremal polynomials for measures
with compact support. However, (\ref{asi-lag-2}) and
(\ref{asi-her-2}) allow us to give the
 outer strong asymptotics of the Laguerre--Sobolev and
 Hermite--Sobolev orthogonal polynomials in terms of the classical
outer strong asymptotics of Laguerre and Hermite polynomials,
respectively. Indeed, there are no differences between case I and
case II in the relative asymptotics if we compare the Sobolev
orthogonal polynomials with an adequate family of polynomials
defined as:
$$R_n(x)=\lim_{\lambda \to \infty}Q_n(x).$$
These polynomials have been known since the pioneering work
\cite{iserles} (see also for more details \cite{margarita}) and
they satisfy $R_n'(x)=k_n\, P_{n-1}(x; \mu_1)$ where $P_{n-1}(x;
\mu_1)$ are the  polynomials orthogonal  with respect  to the
measure $\mu_1$ and $k_n$ depends on the normalization. Now, we
can simplify and list some results about relative asymptotics of
Sobolev orthogonal polynomials:
\begin{proposition} \hspace{3cm} \label{pro-imp}
\begin{enumerate}
\item[(a)] Relative outer asymptotics. $$\lim_{n\to \infty}
\frac{Q_n^{(Li)}(x)}{R_n^{(Li,\infty)}(x)}=\frac{1}{1-\ell}\, ,
\quad i=1,2$$ uniformly on compact subsets of $
\mathbb{C}\setminus \supp(\mu_1), $ where $R_n^{(Li,\infty)}(x)=\\
\lim_{\lambda \to \infty}Q_n^{(Li)}(x), $ $i=1,2,$ and $\ell$ is
given by (\ref{ele}). \item[(b)] Relative outer asymptotics.
$$\lim_{n\to \infty}
\frac{Q_n^{(Hi)}(x)}{R_n^{(Hi,\infty)}(x)}=\theta (\lambda)\,
,\quad i=1,2$$ uniformly on compact subsets of $
\mathbb{C}\setminus \mathbb{R}, $ where
$R_n^{(Hi,\infty)}(x)=\lim_{\lambda \to \infty}Q_n^{(Hi)}(x), $
$i=1,2,$   $\varphi$ and $\theta (\lambda)$ are given by
(\ref{zhu}) and (\ref{theta}), respectively.
\end{enumerate}
\end{proposition}
\textbf{Proof.}  For case I in Laguerre as well as in Hermite
weights, the proof is straightforward because
$R_n^{(L1,\infty)}=L_n^{(\alpha-1)}(x)$ and
$R_n^{(H1,\infty)}(x)=H_n(x).$ On the other hand, for case II  we
have (see  \cite{mei-per-pin} and \cite{alf-mb-rez})
 $$R_n^{(L2,\infty)}(x)=L_n^{(\alpha)}(x)-c_{n-1}\,L_{n-1}^{(\alpha)}(x)=L_n^{(\alpha-1)}(x)-d_{n-1}
  L_{n-1}^{(\alpha)}(x),$$
 with $(R_n^{(L2,\infty)})'(x)=-T_{n-1}(x)=-\left(L_{n-1}^{(\alpha+1)}(x)-
 c_{n-2}\,L_{n-2}^{(\alpha+1)}(x)\right), $ being $T_n$ the polynomials defined in page \pageref{tes}. Since
 $$\frac{Q_n^{(L2)}(x)}{R_n^{(L2,\infty)}(x)}=
 \frac{Q_n^{(L2)}(x)}{L_n^{(\alpha-1)}(x)}\,
 \frac{\sqrt{n}\,L_n^{(\alpha-1)}(x)}{L_n^{(\alpha)}(x)}\,
 \frac{L_n^{(\alpha)}(x)}{T_n(x)}\,
 \frac{T_n(x)}{\sqrt{n}\,R_n^{(L2,\infty)}(x)}$$
 it only remains to use (\ref{asi-lag-2}),
 \cite[Th.4.5, Lem.4.4]{mei-per-pin} and the very well--known relative asymptotics
 for Laguerre polynomials of the different degrees which  is
 deduced from Perron's Theorem (see \cite[p.199]{sz}). Notice that in case I we have
 $\supp(\mu_1)=[0,\infty),$ but in case II associated with
Laguerre weights, $\supp(\mu_1)=[0,\infty)\cup \{-a\}. $ This
ensures that there is not a division by zero when we consider
$\lim_{n\to \infty}L_n^{(\alpha)}(x)/T_n(x)$ (see \cite[Remark
4.6]{mei-per-pin}). On the other hand, the zeros of
$R_n^{(L2,\infty)}(x)$ accumulate on $[0,\infty)$ when $n\to
\infty.$

(b) follows from  Theorem 2  and Proposition 5 in \cite{mb} and
taking into account that $R_n^{(H2,\infty)}(x)$ has not complex
zeros when $n\to \infty $ according to Theorem 6.3 in
\cite{bru-gro-mei}. $\Box$

\section{The Freud weights}

The first example of Sobolev inner product associated with
measures of unbounded support out of the coherent or symmetric
coherent pairs was given in \cite{cac-mar-mb} within  the
framework of the characterization of semiclassical measures. There
the asymptotic properties of the Sobolev polynomials orthogonal
with respect to the inner product,
\begin{equation} \label{pr-freud-sob} (p,q)=\int_{\mathbb{R}}
p(x)q(x)\,e^{-x^4}dx+\lambda\, \int_{\mathbb{R}}
p'(x)q'(x)\,e^{-x^4}dx, \quad \lambda >0,\end{equation} were
deduced. Recently, in \cite{ger-lub-mar} a more general inner
product was considered
\begin{equation} \label{freud-gen}
\left( f,g\right)
=\int_{\mathbb{R}}f(x)g(x)(\psi(x)W(x))^{2}dx+\lambda
\int_{\Bbb{R}}f^{\prime }(x)g^{\prime }(x)W^{2}(x)dx,
\end{equation}
where $\lambda >0,$ $\psi \in L_{\infty}(\mathbb{R})$ is positive
on a set of positive measure, and $W\left( x\right) =\exp \left(
-Q(x)\right) $ is an exponential weight with some special
assumptions over $Q(x)$. The techniques used in this paper are new
due to the lack of an algebraic relation between these Sobolev
orthogonal polynomials and standard polynomials associated with
the weight function $W^2$.

We denote by $p_n=\gamma_n\, x^n+\ldots$ the orthonormal
polynomials associated with $W^2$ and by $q_n=\kappa_n\,
x^n+\ldots$ the Sobolev polynomials orthonormal with respect to
the inner product (\ref{freud-gen}). We also use the
Mhaskar--Rakhmanov-Saff numbers (see, for example,
\cite{mha-saff}), $a_n, \,  n\ge 1,$ which are the positive roots
of the
equations $$ n=\frac{2}{\pi }\int_{0}^{1}a_{n}tQ^{\prime }\left( a_{n}t\right) \frac{dt}{%
\sqrt{1-t^{2}}} $$ where $Q(x)$ is a convex and even function
defined on $\mathbb{R}.$ Notice that in our case $\lim_{n\to
\infty}a_n/n=0\, .$ Furthermore, we need the Szeg\H{o} function
(see \cite{sz}) associated with  a measurable function
$f:\left[-\pi, \pi \right]\to [0,\infty)$ satisfying the condition
$\int_{-\pi }^{\pi }\log f\left( \theta \right) d\theta
>-\infty ,$ that is,
$$
D\left( f;z\right) =\exp \left( \frac{1}{4\pi }\int_{-\pi }^{\pi
}\log f\left( e^{it}\right) \frac{e^{it}+z}{e^{it}-z}dt\right)
,\quad \left| z\right| <1,$$ with $$ \left| D\left( f;e^{i\theta
}\right) \right| ^{2}=f\left( \theta \right) \text{ a.e. }\theta
\in \left[ -\pi ,\pi \right].  $$   The asymptotics of the
polynomials $q_n$ and their derivatives are stated in the next:

\begin{theorem}[Th. 1.1 in \cite{ger-lub-mar}] \label{teo-freud}
\item
 Let $Q:\mathbb{R}\rightarrow \mathbb{R}$  be an even and continuous
function in $\mathbb{R}.$ Assume that $Q^{\prime \prime }$ is
continuous in $\left( 0,\infty \right), $ and $Q^{\prime }>0$ in
$\left( 0,\infty \right). $ Assume, furthermore, that
for some $\alpha ,\beta >0,$%
\begin{equation}
\alpha \leq \frac{xQ^{\prime \prime }\left( x\right) }{Q^{\prime
}\left( x\right) }\leq \beta,\quad x\in \left( 0,\infty \right) .
\end{equation}
 Then, using the previous notation, for $%
n\rightarrow \infty ,$

\begin{itemize}
\item
\begin{equation} \label{siete}
\left\| \left( q_{n}^{\prime }-\frac{1}{\sqrt{\lambda
}}p_{n-1}\right) W\right\| _{L_{2}\left( \Bbb{R}\right) }=O\left(
\frac{a_{n}}{n}\right) =o\left( 1\right)
\end{equation}
and
\begin{equation} \label{ocho}
\left\| \left( 1+\left| Q^{\prime }\right| \right) \left( q_{n}-\frac{1}{%
\sqrt{\lambda }}\int_{0}^{x}p_{n-1}\right) W\right\| _{L_{2}\left( \Bbb{R}%
\right) }=O\left( \sqrt{\frac{a_{n}}{n}}\right).
\end{equation}
\item
\begin{equation} \label{nueve}
\left\| \left( q_{n}^{\prime }-\frac{1}{\sqrt{\lambda
}}p_{n-1}\right)
W\right\| _{L_{\infty }\left( \Bbb{R}\right) }=O\left( \sqrt{\frac{a_{n}}{n}%
}\right).
\end{equation}
and
\begin{equation} \label{diez}
\left\| \left( 1+\left| Q^{\prime }\right| \right) \left( q_{n}-\frac{1}{%
\sqrt{\lambda }}\int_{0}^{x}p_{n-1}\right) W\right\| _{L_{\infty
}\left( \Bbb{R}\right) }=O\left( \sqrt{\frac{a_{n}}{n}}\right).
\end{equation}
\item Let
$$
W_{n}\left( \theta \right) =W\left( a_{n}\cos \theta \right),\quad
\theta \in \left[ -\pi ,\pi \right] .
$$
Uniformly in closed subsets of $\Bbb{C}\backslash \left[
-1,1\right] ,$%
\begin{equation} \label{doce}
\left| \frac{\left( q_{n}^{\prime }-\frac{1}{\sqrt{\lambda }}%
p_{n-1}\right) \left( a_{n}z\right) }{\varphi \left( z\right)
^{n-1}D^{-2}\left( W_{n};1/\varphi \left( z\right) \right)
}\right| =O\left( \frac{\sqrt{a_{n}}}{n}\right).
\end{equation}
Uniformly in closed subsets of $\Bbb{C}\backslash \left[
-1,1\right] ,$%
\begin{equation} \label{trece}
\left| \frac{ q_{n}\left( a_{n}z\right) -\frac{1}{\sqrt{\lambda }}%
\int_{0}^{a_{n}z}p_{n-1} }{\varphi \left( z\right)
^{n}D^{-2}\left( W_{n};1/\varphi \left( z\right) \right) }\right|
=O\left( \frac{a_{n}^{3/2}}{n^2}\right),
\end{equation}
where $\varphi$ is given by (\ref{zhu}).
\end{itemize}
\end{theorem}

\noindent \textbf{Outlines of the proof.}  To establish these
results we need to introduce some concepts. Let $f$ be an
absolutely continuous function, then
$$ E_{n}\left[ f;W\right] =\inf \left\{
\parallel (f-P)W\parallel _{L_{2}\left( \mathbb{R}\right) }:\deg \left(
P\right) \leq n\right\}
$$
is the weighted $L_{2}$ error of approximation. Let
$$
I\left[ R\right] \left( x\right) =W\left( x\right)
^{-2}\int_{-\infty}^{x}R\left( t\right) W^{2}\left( t\right) \psi
^{2}\left( t\right) dt,\quad x\in \mathbb{R} ,
$$
be the linear operator defined on suitably restricted classes of
functions $R$.

We begin giving the main ideas of the proof of (\ref{siete}).
\begin{itemize}
\item In the framework of more general weights (see \cite[proof of
Th. 1.3]{ger-lub-mar}) an upper bound for the ratio
$\gamma_{n-1}/\gamma_n$ depending of $n$ is proved, that is,
\begin{equation} \label{des-gammas}
\frac{\gamma_{n-1}}{\gamma_n}\leq 4a_n\, .
\end{equation}
\item Using the Jackson inequality for exponential weights (see,
for example, \cite{dit-lub})
$$E_n[f;W]\le C\frac{a_n}{n}\, ||f^{\prime}
W||_{L_2(\mathbb{R})}, $$ the statement of Theorem 1.5 in
\cite{ger-lub-mar} asserts that
$$
\left| \left( \frac{\gamma _{n-1} }{n\kappa _{n}}%
\right) ^{2}-\lambda \right|
\leq \left( \parallel \psi \parallel _{L_{\infty }\left( \mathbb{R}\right) }\frac{%
\gamma _{n-1} }{n\gamma _{n} }\right)
^{2}+C_{1}\left(\frac{a_{n-2}}{n-2}\right)^{2}.
$$
Now, using (\ref{des-gammas}) and the properties of $\psi$ we
obtain
\begin{equation} \label{88}
\left| \left( \frac{\gamma _{n-1} }{n\kappa _{n}}%
\right) ^{2}-\lambda \right|=O\left( \frac{a_n}{n} \right)^2.
\end{equation}
\item From the Sobolev orthogonality of the polynomials $q_n$  we
get
\begin{equation} \label{79}
\int_{\mathbb{R}}\left( q_{n}^{\prime }-n\frac{\kappa _{n}}{\gamma
_{n-1} }p_{n-1}\right) ^{2}W^{2} \leq \frac{1}{\lambda }-\left( n%
\frac{\kappa _{n}}{\gamma _{n-1}}\right) ^{2}
\end{equation}
\item Finally, taking into account
\begin{align*} \int_{\mathbb{R}}\left( q_{n}^{\prime
}-\frac{1}{\sqrt{\lambda }}p_{n-1}\right) ^{2}W^{2}
=&\int_{\mathbb{R}}\left( q_{n}^{\prime }-n\frac{\kappa
_{n}}{\gamma
_{n-1} }p_{n-1}\right) ^{2}W^{2}\\&+\int_{\mathbb{R}}\left( n\frac{%
\kappa _{n}}{\gamma _{n-1} }-\frac{1}{\sqrt{\lambda }}%
\right) ^{2}p_{n-1}^{2}W^{2},
\end{align*}

it only remains to use (\ref{88}--\ref{79}) in order to obtain
(\ref{siete}).
\end{itemize}

(\ref{ocho}) follows from four results: (\ref{siete}),
\cite[Lem.3.4.(c)]{ger-lub-mar}, that is, \begin{equation}
\label{5.2}
\parallel \left( 1+\left| Q^{\prime }\right| \right) \left(
q_{n}-q_{n}\left( 0\right) -\frac{1}{\sqrt{\lambda }}\int_{0}^{x}p_{n-1}%
\right) W\parallel _{L_{2}\left( \Bbb{R}\right) }=O\left( \frac{a_{n}}{n}%
\right),
\end{equation}
 the bound \begin{equation} \label{cota-cero} |q_n(0)|\le C\sqrt{\frac{a_n}{n}} \end{equation} given in
 \cite[Lem.5.2]{ger-lub-mar}, and the fact that $\parallel \left( 1+\left| Q^{\prime }\right| \right) W\parallel
_{L_{2}\left( \Bbb{R}\right) }$ is finite.

To prove (\ref{nueve}) we only need (\ref{siete}) and the
Nikolskii inequality $\parallel
PW\parallel_{L_{\infty}(\mathbb{R})}\le C
\left(\frac{n}{a_n}\right)^{1/2} \parallel
PW\parallel_{L_{2}(\mathbb{R})}.$ The proof of (\ref{diez})
follows the proof of (\ref{ocho}) but now we use (\ref{nueve}) and
$$
\parallel \left( 1+\left| Q^{\prime }\right| \right) \left(
q_{n}-q_{n}\left( 0\right) -\frac{1}{\sqrt{\lambda }}\int_{0}^{x}p_{n-1}%
\right) W\parallel _{L_{\infty}\left( \Bbb{R}\right) }=O\left( \sqrt{\frac{a_{n}}{n}}%
\right). $$

Finally, (\ref{doce}) follows from (\ref{siete}) and applying
Lemma 14.6 in \cite[p.395]{lev-lub} with $p=2, m=n-1,$ and
$w=W_n.$ To establish (\ref{trece}) we need this Lemma 14.6,
(\ref{5.2}), (\ref{cota-cero}), and some other technical results.

From this theorem and very well--known results for orthogonal
polynomials associated with exponential weights, see for example
\cite{lev-lub}, we deduce:

\begin{corollary}[Cor. 1.2 in \cite{ger-lub-mar}] \hspace{3cm} \begin{enumerate}
\item[(a)] When $n$ tends to infinity
$$
\kappa _{n}=\frac{1}{n}\frac{1}{\sqrt{2\pi \lambda }}\left( \frac{a_{n}}{2}%
\right) ^{-n+\frac{1}{2}}\exp \left( \frac{2}{\pi
}\int_{0}^{1}\frac{Q\left( a_{n}s\right)
}{\sqrt{1-s^{2}}}ds\right) \left( 1+o\left( 1\right) \right) .
$$
\item[(b)] For $n\rightarrow \infty ,$%
\begin{align*}
\int_{-1}^{1}\left| \sqrt{\lambda a_{n}}\right. &q_{n}^{\prime
}\left(
a_{n}x\right) W(a_{n}x)-\frac{\sqrt{2/\pi }}{(1-x^{2})^{1/4}} \\
& \left.\times \cos \left[ \left(n-\frac{1}{2}\right)\arccos
x+2\Gamma (W_{n};\arccos x)-\frac{\pi }{4}\right] \right|
^{2}dx=o\left( 1\right) . \end{align*}
\item[(c)] Uniformly for $z$ in closed subsets of $\Bbb{C}%
\backslash [-1,1],$ we have, as $n\rightarrow \infty ,$%
\begin{align*}
\sqrt{\lambda a_{n}}q_{n}^{\prime }(a_{n}z)/&\left\{ \varphi
(z)^{n-1}D^{-2}\left(W_{n};\frac{1}{\varphi (z)}\right)(1-\varphi (z)^{-2})^{-1/2}\right\} \\&=%
\frac{1}{\sqrt{\pi }}(1+o\left( 1\right) ). \end{align*}
\item[(d)] There exists $\eta >0$ such that as $n\rightarrow
\infty ,$ we have, uniformly for $\left| x\right| \leq 1-n^{-\eta
},x=\cos \theta ,$%
\begin{align*}
\sqrt{\lambda a_{n}}q_{n}^{\prime
}(a_{n}x)W(a_{n}x)(1-x^{2})^{1/4}& =\sqrt{\frac{2}{\pi }}\cos
\left( \left(n-\frac{1}{2}\right)\theta +2\Gamma (W_{n};\theta
)-\frac{\pi }{4}\right)\\ & +O(n^{-\eta }) \end{align*}
\end{enumerate}
\end{corollary}

\noindent \textbf{Remarks.}

\begin{itemize}
\item  Theorem \ref{teo-freud} can be formulated for more general
weights, see Theorem 1.3 in \cite{ger-lub-mar}, even for weights
such that the asymptotics for the corresponding orthogonal
polynomials is unknown, and therefore the previous Corollary
cannot be established.

\item The main key to obtain (\ref{siete}) and, therefore some of
the other asymptotic properties of the polynomials $q_n\, ,$ is
\cite[formula (34) in Th.1.5]{ger-lub-mar}. To obtain this formula
we need a Jackson inequality for exponential weights and bounds of
the sequence $\left( \gamma_{n-1}/(n\kappa_n) \right)^2$ given in
Theorem 1.4 of \cite{ger-lub-mar}, i.e,
\begin{equation} \label{des-lideres}
\lambda+\left( \frac{\gamma_{n-1}}{n\gamma_n^*} \right)^2\le
\left( \frac{\gamma_{n-1}}{n\kappa_n} \right)^2\le \lambda+\left(
\frac{\gamma_{n-1}}{n\gamma_n^*} \right)^2+\sup E_{n-2}^2[I[R];W],
\end{equation}
where $\gamma_n^*$ is the leading coefficient of the orthonormal
polynomials $\pi_n$ for the weight function $(\psi\, W)^2,$ and
where the $\sup $ is taken over all polynomials $R$ of degree
$\leq n-1$ satisfying both
$$
\parallel RW\psi \parallel _{L_{2}\left( \mathbb{R}\right) }=1 \quad
\textrm{and } \int_{\mathbb{R}}R\left( W\psi \right) ^{2}=0.
$$
Theorem 1.4 in \cite{ger-lub-mar} is essential to establish
Theorem 1.1 as well as Theorem 1.3 for more general weights in
\cite{ger-lub-mar}. So, it is important to point out that the
inequalities given by (\ref{des-lideres}) for the leading
coefficients are equivalent to the ones for the square of the
norms used to establish the asymptotics of the Sobolev orthogonal
polynomials in the framework of the coherence (see, for example,
\cite[Th.2]{and-mb-per-pin}, \cite[p.52]{alf-and-rez},
\cite[Prop.3]{mar-mb} or \cite[proof of Lemma
2.6]{alf-mb-per-pin-rez}). Furthermore, these inequalities for the
square of the norms can be obtained when we have a finite
algebraic relation between Sobolev polynomials and the standard
ones (one example is given in \cite{cac-mar-mb}).

On the other hand, the lower bound for $\left(
\gamma_{n-1}/(n\kappa_n) \right)^2$ given in (\ref{des-lideres})
is a ``universal" property of the Sobolev orthogonality, that is,
if we consider the monic polynomials $\Pi_n=\pi_n/\gamma_n^*\, ,$
$ P_n=p_n/\gamma_n\, ,$ and $Q_n=q_n/\kappa_n\,$ we have,
\begin{align*}
(Q_n,Q_n)&=\int Q_n^2(x)\, (\psi
W)^2(x)dx+\lambda \int (Q_n^{\prime})^2(x)\,  W^2(x)dx \\
&\ge \int \Pi_n^2(x)\, (\psi W)^2(x)dx+\lambda\, n^2 \int
P_{n-1}^2(x)\,  W^2(x)dx.
\end{align*}
where  the extremal property of the square of the norms of $\Pi_n$
and $P_n$ have been used. Now, this relation can be expressed in
terms of the leading coefficients as $$\frac{1}{\kappa_n^2}\ge
\left(\frac{1}{\gamma_n^*}\right)^2+\lambda\,
\left(\frac{n}{\gamma_{n-1}}\right)^2, $$ from where the
inequality on the left of (\ref{des-lideres}) follows. However,
until the work \cite{ger-lub-mar}, an upper bound for $\left(
\gamma_{n-1}/(n\kappa_n) \right)^2$ had been obtained only
assuming the existence of a finite algebraic relation between
Sobolev and standard polynomials. Therefore, the upper bound given
in (\ref{des-lideres}) is an important advance in the study of
Sobolev orthogonal polynomials.
\end{itemize}

There are not a lot of results about the zeros of these Sobolev
polynomials for Freud weights. For the polynomials orthogonal with
respect to the particular inner product (\ref{pr-freud-sob}) it
has been proved that the zeros, $q_{n,k}^{\lambda}\, ,$ of the
corresponding Sobolev orthogonal polynomials are real, simple and
interlace with those of the standard polynomials orthogonal for
the weight function $e^{-x^4}$ on $\mathbb{R}$  (see
\cite{mb-pr}). Are the zeros of Sobolev polynomials for
exponential weights like those considered in \cite{ger-lub-mar}
real and simple? This question remains open.

On the other hand, in a very recent work \cite{lop-mar-pij} the
conditions on the weights in \cite{ger-lub-mar} have been relaxed
 aiming for weak asymptotic instead of strong
asymptotic results. In fact, for $p\in[1,\infty)\,$  the authors
consider a family of weights $\{w_i\}_{i=0}^{m}$ and  the Sobolev
$p$--norm \begin{equation} \label{p-norm}
||q||_S=\left(\sum_{k=0}^{m}
\int_\mathbb{R}|q^{(k)}(x)w_k(x)|^pdx\right)^{1/p}=\left(
\sum_{k=0}^{m} ||q^{(k)}w_k||_{L^p(\mathbb{R})}^{p}
\right)^{1/p}\, , q\in \mathbb{P}\, ,\end{equation} where
$\mathbb{P}$ is the linear space  of polynomials with real
coefficients. Then, it is said that $Q_n$ is a $n$--th extremal
monic polynomial with respect to (\ref{p-norm}) if
$$||Q_n||_S=\min \{||q||_S :q(x)=x^n+\ldots \}\, .$$ The
existence of such a type of polynomials is easy to prove when
$n\in \mathbb{Z}_+\, .$ Next, we use the following definition:
\begin{definition} Let $w$  be a positive continuous function on
$\mathbb{R}$.
\begin{itemize}
\item $w\in W(\alpha, \tau)$ with $ \alpha, \tau >0,$ if
$$\lim_{|x|\to \infty} \frac{-\log w(x)}{\tau\, |x|^{\alpha}}=1.$$
\item $w\in W(\alpha)$ with $ \alpha>0,$ if
$$\lim_{|x|\to \infty} \frac{\log\log\frac{1}{
w(x)}}{\log|x|}=\alpha.$$
\end{itemize}
\end{definition}
Notice that $W(\alpha, \tau)\subset W(\alpha), $ for all $\tau>0.$
In particular, the Freud--type weights $e^{-\tau |x|^{\alpha}}\in
W(\alpha, \tau)\, .$

Let $\overline{k}\in \{0,1,\ldots, m\}$ be the smallest index such
that either $$\alpha_{\overline{k}}<\min_{\substack{ 0\leq k\leq m
\\ k\neq \overline{k}}} \alpha_k\, ,\quad \mathrm{or}
\quad
\alpha_{\overline{k}}=\min_{\substack{0\leq k\leq m \\
k\neq \overline{k} }} \alpha_k \, \quad \mathrm{and} \, \quad
\tau_{\overline{k}}=\min \{\tau_k: \alpha_k=\alpha_{\overline{k}}
\}, $$ and let $\widetilde{k}\in \{0,1,\ldots, m\}$ be the first
index such that $$\alpha_{\widetilde{k}}=\min\{\alpha_k : 0\le
k\le m\}.$$ We also use the notation
$\overline{\alpha}=\alpha_{\overline{k}}\, ,$
$\overline{\tau}=\tau_{\overline{k}}\, ,$ and
$\widetilde{\alpha}=\alpha_{\widetilde{k}}\, .$ In
\cite{lop-mar-pij}, the $n$--th root and logarithmic  asymptotics
of $Q_n$ are established:
\begin{theorem}[Th. 1.1 and 1.2 in \cite{lop-mar-pij}]
\item
\begin{enumerate}
\item[(a)] Let $Q_n$ be the $n$--th Sobolev monic extremal
polynomials relative to the norm (\ref{p-norm}) where $w_k \in
W(\alpha_k, \tau_k),$ $\, 0\le k\le m.$ Then
$$\lim_{n\to \infty}
n^{-1/\overline{\alpha}}||Q_n||_S^{1/n}=\frac{1}{2} \, \left(
\frac{\gamma_{\overline{\alpha}}}{\overline{\tau}\, e
}\right)^{1/\overline{\alpha}}, $$ and, for all $j\ge
\overline{k}$ \begin{equation}  \label{eq:j-der} \lim_{n\to
\infty} n^{-1/\overline{\alpha}}||Q_n^{(j)}e^{-\overline{\tau}
|x|^{\overline{\alpha}}}||_{L_{\infty
(\mathbb{R})}}^{1/n}=\frac{1}{2} \, \left(
\frac{\gamma_{\overline{\alpha}}}{\overline{\tau}\, e
}\right)^{1/\overline{\alpha}}, \end{equation} where $$
\gamma_{\overline{\alpha}}=\frac{\Gamma\left(\overline{\alpha}/2\right)\,
\Gamma (1/2)}{2\, \Gamma\left( (1+\overline{\alpha})/2\right)} \,
.$$

\item[(b)] Let $Q_n$ be the $n$--th Sobolev monic extremal
polynomials relative to the norm (\ref{p-norm}) where $w_k \in
W(\alpha_k),$ $\, 0\le k\le m.$ Then
$$\lim_{n\to \infty}
||Q_n||_S^{1/(n\log n)}=e^{1/\widetilde{\alpha}}\, .$$

\end{enumerate}
\end{theorem}

\noindent \textbf{Remarks.}

\begin{itemize}
\item The asymptotic contracted limit distribution of the zeros of
$Q_n^{(j)}$ and the weak limit of the corresponding contracted
extremal polynomials is obtained  in \cite[Cor.3.1 and
3.2]{lop-mar-pij} as a consequence of (\ref{eq:j-der}).

\item Some of the techniques used to prove these results are
similar to those used in \cite{lop-pij} and \cite{lop-pij-per} for
the bounded case.

\end{itemize}

\section{Some open problems and new directions}

In the previous sections some unsolved questions  and conjectures
have been posed in their natural context. Next, we wish to give
some guidelines which we consider as new directions in the study
of the analytical properties of Sobolev orthogonal polynomials
associated with measures with unbounded support.

\bigskip

\noindent \textbf{1.}  \cite{ger-lub-mar} permits the study of
analytic properties of the polynomials orthogonal with respect to
the inner product
$$(p,q)=\int p(x)q(x)d\mu_0+\lambda \int_{0}^{\infty}
p^{\prime}(x)q^{\prime}(x)x^{\alpha}e^{-x}dx, \quad \alpha>-1,$$
where the measure $\mu_0$ is in some sense ``admissible", for
example, $d\mu_0(x)=\psi(x)\, x^{\alpha}e^{-x}dx$ where $\psi$
satisfies suitable conditions. Under these assumptions it is
possible to establish an analog  of Theorem 1.4 in
\cite{ger-lub-mar}.

\bigskip

\noindent \textbf{2.} In this paper we have dealt with a Sobolev
inner product involving the derivative operator. However, another
possibility is to consider the inner product
\begin{equation} \label{int-3}
(p, q) = \int_{ \mathbb{R}} p(x)q(x)\, d\mu_0 (x)+\lambda\,\int_{
\mathbb{R}} \left( \Delta p\right)(x) \left( \Delta q \right)(x)\,
d\mu_1 (x)\, ,
\end{equation}
where $\mu_0$, $\mu_1$ are measures supported on a countable set
and $(\Delta p)(c)=p(c+1)-p(c)$ denotes the forward difference
operator. The goal of this type of inner product is to show that a
family of continuous Sobolev orthogonal polynomials can be seen as
a limit of the $\Delta$--Sobolev polynomials orthogonal with
respect to the inner product (\ref{int-3}). Notice that in
\cite{are-god-mar-mb} an inner product like (\ref{int-3}) with
$\mu_1$ as the measure corresponding to the Meixner weight
function was considered in the framework of $\Delta$--coherence.
More precisely, according to the classification of
$\Delta$--coherence pairs (see \cite{are-god-mar}) if $\mu_1$ is a
classical discrete measure (in our case, Meixner) then $\mu_0$ is
a polynomial modification of degree one of the measure $\mu_1.$

We denote by $M_{n}^{(\gamma, \mu)}(x)$ the monic Meixner
polynomials orthogonal with respect to
$$(p,q)=\sum_{x=0}^\infty p(x)q(x) \frac{\mu^x\,
\Gamma(x+\gamma)\,(1-\mu)^\gamma}{\Gamma(\gamma) \Gamma(x+1)} \, ,
\quad \gamma>0\, , \quad 0<\mu<1, $$ and by
$Q_n(x;\gamma,\mu;\lambda,K,a)$ the $\Delta$--Meixner--Sobolev
polynomial orthogonal with respect to (\ref{int-3}) where
$(\mu_0,\mu_1)$ is a $\Delta$--coherent pair of measures of case I
(Meixner), that is, $\mu_1$ is the Meixner measure. Notice that
the parameters $a, K$ appear in the expression of the measure
$\mu_0$ (see \cite[Prop. 2.10]{are-god-mar-mb}). In
\cite{are-god-mar-mb} the outer asymptotics of these polynomials
$Q_n(x;\gamma,\mu;\lambda,K,a)$ was studied. Furthermore, some
algebraic and analytic results for coherent pairs associated with
Laguerre weights were recovered in Section 4 of that work. In
particular,  denoting by $Q_n(x;\alpha;\lambda,M,a):=(-1)^n\,
n!\,Q_n^{(L1)}(x)$ being $Q_n^{(L1)}(x)$ the Laguerre--Sobolev
orthogonal polynomials defined in Theorem \ref{ceros-lag-1}, then
(see \cite[Prop. 4.4]{are-god-mar-mb})
\begin{equation}\label{E:LMLS}
\lim_{\mu \uparrow 1} (1-\mu)^n Q_n\left( \frac{x}
{1-\mu};\alpha+1,\mu;\frac{\lambda}{(1-\mu)^2},M,\frac{a}{1-\mu}
\right)=Q_n(x;\alpha;\lambda,M,a)\, .
\end{equation}
 (\ref{E:LMLS}) is a limit relation between
 $\Delta$--Meixner--Sobolev orthogonal polynomials and the
 Laguerre--Sobolev ones which is similar to the one existing
 between Meixner and Laguerre polynomials in the Askey scheme (see,
 for example, \cite{and-ask}). A natural consequence is the
 possibility of establishing  a Sobolev--Askey scheme.

\bigskip

\noindent \textbf{3.} As we have commented previously, the second
measure $\mu_1$ plays the main role in the Sobolev orthogonality.
In order to equalize the role of both measures, the balanced
Sobolev orthogonal polynomials were introduced in
\cite{alf-and-rez} where both measures  have compact support. It
seems natural that we can use  this approach for the study of
Sobolev orthogonal polynomials in the unbounded case, although
this presents some complications. We are working in this
direction, and we have a new result where we show how to balance
in the unbounded case (see \cite{alf-mb-pen-rez}).

\bigskip

\noindent \textbf{4.} We consider the Sobolev--type inner product:
\begin{equation} \label{type-gen}
(p,q)=\int_{I} p(x) q(x)d\mu+\sum_{i=0}^{N}M_i p^{(i)}(c)\,
q^{(i)}(c), \quad M_i\ge 0, \, \, c\in \mathbb{R},
\end{equation}
where $\mu$ is a positive measure with unbounded support
$I\subseteq \mathbb{R}.$ Algebraic results for the orthogonal
polynomials can be deduced regardless  whether the measure $\mu$
is unbounded or not. The study of the asymptotics for the
polynomials orthogonal with respect to (\ref{type-gen}) has been
done especially when $\mu$ is bounded (see the survey
\cite{and-survey}), and recently in \cite{ran-mb} when $\mu$ is
the measure associated with   the Laguerre weight, $N=1,$ and
$c=0.$ Therefore, the study of the asymptotics  and the zero
asymptotics of the polynomials orthogonal with respect to
(\ref{type-gen}) when the measure $\mu$ is more general, for
example associated with  Freud weights, remains open.

\bigskip

\noindent \textbf{5.} Finally, a very interesting problem is to
consider the inner product:
\begin{equation} \label{pr-super-ge}
(p,q)=\int_{I_0} p(x)\, q(x)d\mu_0(x)+\lambda \int_{I_1}
p^{\prime}(x)\, q^{\prime}(x)d\mu_1(x),
\end{equation}
where $\mu_i$ are positive measures supported on $I_i\subseteq
\mathbb{R}, $ $i=1,2,$ respectively, with non--zero absolutely
continuous parts,  and $I_0\cap I_1=\emptyset$ or, at least,
$I_0\supset I_1.$ The main work in this direction is
\cite{gau-kui} where the authors consider measures with compact
support under some general assumptions and they study the
asymptotic distribution of the zeros and critical points of the
polynomials orthogonal with respect to (\ref{pr-super-ge}), and
 they formulated some conjetures. As far as we know,
there are no more works in this direction, and therefore this is
an open problem for both the unbounded and the bounded cases.

\bigskip

\noindent \textbf{Acknowledgements.} The authors thank the referee
by the careful revision of the manuscript and the suggestions for
improve the presentation.

\end{document}